\documentclass[11pt,reqno]{amsart}
\usepackage{amsmath, amssymb, amsthm, microtype,mathtools, bbm}
\usepackage[dvipsnames]{xcolor}
\usepackage{arydshln, comment,dsfont}
\usepackage{verbatim}
\usepackage{graphicx}
\usepackage{yhmath}
\usepackage{geometry}
\usepackage{thmtools}
\usepackage{thm-restate}
\usepackage[hidelinks]{hyperref}
\usepackage[msc-links, initials]{amsrefs}
\DeclareMathOperator{\lcm}{lcm}
\makeatletter
\renewcommand{\pod}[1]{\allowbreak\mathchoice
  {\if@display \mkern 18mu\else \mkern 8mu\fi (#1)}
  {\if@display \mkern 18mu\else \mkern 8mu\fi (#1)}
  {\mkern4mu(#1)}
  {\mkern4mu(#1)}
}
\newcommand{\pdiv}{\mid\!\mid}

\renewcommand{\pmod}[1]{\allowbreak\mkern6mu(\mathrm{mod}\mkern6mu#1)}
\newtheorem{theorem}{Theorem}[section]
\newtheorem{proposition}[theorem]{Proposition}

\newtheorem{lemma}[theorem]{Lemma}
 \newtheorem{fact}[theorem]{Fact}
 \newtheorem{question}[theorem]{Question}

\theoremstyle{definition}
\newtheorem*{definition*}{Definition}
\newtheorem{observation}[theorem]{Observation}

\theoremstyle{remark}
\newtheorem{remark}[theorem]{Remark}
\numberwithin{equation}{section}

\title{Unit fraction representations using shifted Golomb primes}
\author{Rishika Agrawal}
\address{Department of Mathematics, University of Georgia, Athens, GA 30602, USA}
\email{rishika.agrawal@uga.edu}

\begin{document}

\maketitle
\begin{abstract}
We prove that every positive rational number is a sum of distinct unit fractions with denominators of the form $g-h$, where $g$ is a Golomb prime and $h$ is a fixed nonzero integer, as long as the necessary congruence conditions are satisfied. Golomb primes, introduced by Golomb in 1955, are a natural analogue of the
primes.
Erd\H{o}s showed that the relative density of  Golomb primes among the
primes up to $N$ is about $1/\log\log N$, so they form a
sparse subset of the primes. Our result is an analogue of Bloom's
result on unit fraction representations using shifted primes, for this sparse subset of the primes. A key ingredient in our proof is a Siegel--Walfisz-type theorem for Golomb primes in arithmetic progressions, which we also prove here.

\end{abstract}

\section{Introduction}

A unit fraction is a fraction of the form $1/n$, where $n$ is a positive integer. A \emph{unit fraction representation} of a rational number is a way of writing it as a sum of finitely many distinct unit fractions. The greedy algorithm gives an elementary proof that every positive rational number has such a representation (see, for instance,~\cite{Sigler2002}). Questions about unit fractions are among the oldest in number theory, and many of them are still open (see the surveys~\cite{Graham2013Egyptian, BloomElsholtz2022}).

In this paper we study unit fraction representations in which every denominator lies in a given set $A\subseteq\mathbb{N}$. Many such results are
known, see for example~\cite{Breusch1954,Stewart1954,AlbadaLint1963,Graham1964,ErdosGraham1980,ButlerErdosGraham2015,Eppstein2021, vanDoornKovac2026}.
A well-known problem of this kind is the coloring conjecture of Erd\H{o}s and Graham~\cite{ErdosGraham1980}, which says that if the integers greater than $1$ are colored with finitely many colors, then there exists a monochromatic finite subset $B$ of integers such that $1=\sum_{n\in B}1/n$. Croot~\cite{Croot2003} proved this conjecture, and Bloom~\cite{Bloom2025density} later proved a stronger density version. More precisely, Bloom showed that every set $A\subseteq[1,N]$ with \[\sum_{n\in A} \frac{1}{n}\gg\dfrac{{\log \log \log N}}{\log \log N}\log N\] contains a subset $B$ with $\sum_{n\in B}1/n=1$. Liu and Sawhney~\cite{LS24} then relaxed the requirement to $\sum_{{n\in A}}1/n\geq(\log N)^{4/5+o(1)} $. A natural question studied by Graham~\cite{Graham1964} (see also~\cite{Wilf1961}) and raised again by Bloom~\cite{Bloom2025shifted}, is the following.

\begin{question}\label{q:bloom}
For which sets $A\subseteq\mathbb{N}$ can every positive rational number (subject to the necessary congruence conditions) be written as $\sum_{n\in B}1/n$ for some finite set $B\subseteq A$?
\end{question}

In the same paper~\cite{Bloom2025shifted}, Bloom showed that for every fixed nonzero integer $h$, the shifted primes $\{p-h: p \text{ prime},\ p>h\}$ have the property in Question \ref{q:bloom}, proving a conjecture of Sun~\cite[Conjecture~4.1]{Sun2017}. However, the primes themselves do not have this property; for example, $1$ cannot be written as a sum of reciprocals of distinct primes. Pomerance ~\cite{Croot2000,Bloom2025density} constructed another set without this property, by taking the set of all $n$ whose largest prime factor is at least $cn/\log n$, for a suitable constant $c>0$. This set is fairly large, in the sense that the sum of the reciprocals of its elements $n\leq N$ is of order $(\log\log N)^2$, whereas for the primes this sum is about $\log\log N$.

In the other direction, one can ask how thin a set can be and still have this property. Our main result is that the set of shifted Golomb primes, defined below, has this property. The Golomb primes form a very sparse subset of primes, Erd\H{o}s \cite{Erdos1961} showed that the sum of their reciprocals up to $N$ is only about $\log\log\log N$.

To motivate the definition of Golomb primes, it helps to first describe the primes as the output of a sieve. Set $p_1=2$. Once $p_1,\dots,p_k$ are defined, let $p_{k+1}$ be the least natural number greater than $p_k$ such that $p_{k+1}\not\equiv 0\pmod{p_i}$ for all $1\leq i\leq k$. This is the sieve of Eratosthenes. Golomb~\cite{Golomb1955} used a similar sieve on the primes, removing the residue class $1$ in place of the class $0$. His motivation came from studying integers $n$ with $(n,\phi(n))=1$, called cyclic numbers~\cite{Erdos1948, Pollack2022Cyclic}. These are exactly the $n$ for which every group of order $n$ is cyclic. Every product of distinct Golomb primes is a cyclic number.

\begin{definition*}
The \emph{Golomb primes} are defined inductively as follows. Let $g_1:=3$. Given $g_1,\dots,g_k$, let $g_{k+1}$ be the least prime greater than $g_k$ such that $g_{k+1}\not\equiv 1\pmod{g_i}$ for all $1\leq i\leq k$. We write $\mathcal{G}$ for the set of Golomb primes, so
\[
  \mathcal{G}=\{3,5,17,23,29,53,83,89,113,\dots\}.
\]
\end{definition*}
The Golomb primes are a natural analogue of the primes, but apart from the work of Erd\H{o}s~\cite{Erdos1961}, little is known about them. Let $\mathcal{G}(x)=\mathcal{G}\cap[1,x]$ be the set of Golomb primes up to $x$. Erd\H{o}s~\cite{Erdos1961} proved the following:

\begin{theorem}[Erd\H{o}s~\cite{Erdos1961}]\label{thm:erdos}
We have

\begin{equation}\label{eq:Gx}
    |\mathcal{G}(x)| =(1+o(1)) \frac{x}{\log x\,\log \log x}
\end{equation}
and
\begin{equation} \label{eq:q-1}
    \prod_{q\in\mathcal{G}(x)}\left(1-\frac{1}{q-1}\right)=\frac{1+o(1)}{\log \log{x}}.
\end{equation}

Further, there exists a constant \(c>0\) such that
\begin{equation}\label{eq:1q}
    \sum_{q\in\mathcal{G}(x)}\frac{1}{q} = \log \log \log x-c+o(1).
\end{equation}
 \end{theorem}
   By equation \eqref{eq:Gx}, the Golomb primes up to $x$ make up only about a $1/\log \log x$ proportion of the primes up to $x$, so in particular they have relative density zero in the primes.
    Hence, Bloom's result on shifted primes~\cite{Bloom2025shifted}, which requires positive relative density, does not apply to the Golomb primes. On the other hand, by equation \eqref{eq:1q} the sum of reciprocals of Golomb primes up to $x$ is  about $\log \log \log x$, so it  still diverges. This makes them a natural candidate for a sparse set with the property in Question~\ref{q:bloom}. Our main result shows that they do have this property.
    We say that an integer is \emph{Golomb-prime-free} if no Golomb prime divides it. More precisely, we prove the following.
\begin{restatable}[Main result]{theorem}{mtheoremone}\label{thm:main}
    Let $h$ be a nonzero integer, and let $a,b\in \mathbb{Z}^{+}$ be such that $(a,b)=1$, $(b,h)=1$ and $(b,h-1)$ is Golomb-prime-free. (When $h=1$, this means that no Golomb prime divides $b$.) Then, for every $y\geq|h|$, there is a finite subset $B$ of Golomb primes, all greater than $y$, such that
    \begin{equation}\label{eq:Golrep}
     \frac{a}{b}=\sum_{g\in B}\frac{1}{g-h}.\end{equation}

\end{restatable}
The conditions on $b$ are necessary because of congruence obstructions (see the beginning of Section~\ref{mainthm}).

To prove Theorem~\ref{thm:main}, we need information on the distribution of Golomb primes in arithmetic progressions. We prove the following Siegel--Walfisz-type theorem for Golomb primes, which may be of independent interest.

Let $\mathcal{G}(x;b,h)$ be the set of all Golomb primes $g\leq x$ with $g\equiv h\pmod b$.
\begin{restatable}[Golomb primes in arithmetic progressions]{theorem}{mtheoremtwo}
\label{thm:SW}
    Fix $A>0$. For all $1\leq b\leq(\log x)^A$ and all nonzero $h\in \mathbb{Z}$ such that $(b,h)=1$ and $(b,h-1)$ is Golomb-prime-free, we have
    \[
    |\mathcal{G}(x;b,h)|=(1+o(1)) \frac{x}{\phi (b)\log x\log \log x}\prod_{\substack{q\mid b\\q\in \mathcal{G}}}\left( 1- \frac{1}{q-1} \right)^{-1},
    \]
where the $o(1)$ term depends on $A$.
\end{restatable}

Theorem~\ref{thm:SW} says that the Golomb primes are equidistributed among the allowed residue classes $h\pmod{b}$. The factor $x/(\log x\log\log x)$ is the asymptotic number of Golomb primes up to $x$ (Theorem~\ref{thm:erdos}), and the factor $1/\phi(b)$ is the same as for primes in arithmetic progressions~\cite[Chapter~4]{Davenport2000}. Now suppose that a Golomb prime $q$ divides $b$. Then $q\nmid h-1$, since $(b,h-1)$ is Golomb-prime-free. Every prime $p\equiv h\pmod b$ also satisfies $p\equiv h\pmod q$, and hence $p\not\equiv 1\pmod q$. So for every Golomb prime $q\mid b$, the condition $p\not\equiv 1 \pmod{q}$ already holds for every prime  $p\equiv h\pmod{b}$, and this gives the extra factor $\prod_{q\mid b,\,q\in\mathcal{G}}\big(1-\frac{1}{q-1}\big)^{-1}$.

We call the positive integers $g-h$ with $g\in\mathcal{G}$ the {$h$-shifted Golomb primes}. To prove Theorem~\ref{thm:main}, we also need to know that not too many of the $h$-shifted Golomb primes are non-smooth, that is, have a large prime-power divisor. Let
\[
G_{\rm nsm}(x;b,h;z):=\#\{g\in\mathcal{G}(x;b,h):\ g>h\ \text{and}\ g-h\ \text{has a prime-power divisor greater than}\ z\}.
\]

 \begin{restatable}{theorem}{mtheoremthree}\label{thm:SWnonsmooth}
    Fix $A>0$. For all sufficiently large $x$, all integers $h\not=0$ with $|h|\leq x^{0.01}$  and all $1\leq b\leq (\log x)^A$ such that $(b,h)=1$ and $(b,h-1)$ is Golomb-prime-free, we have \[{G}_{\rm nsm}(x;b,h;x^{0.99})\leq (0.259+o(1)) \frac{x}{\phi(b)\log x\log\log x}\prod_{\substack{q\mid b\\q\in \mathcal{G}}}\left( 1-\frac{1}{q-1} \right)^{-1}, \]
    where the $o(1)$ term depends on $A$.
\end{restatable}
A comparison with Theorem~\ref{thm:SW} shows that at most about $26\%$ of the Golomb primes $g\leq x$ with $g\equiv h\pmod b$ have $g-h$ divisible by a prime power greater than $x^{0.99}$. So a positive proportion of the $h$-shifted Golomb primes remain after we remove the non-smooth ones.
The proofs of Theorems~\ref{thm:SW} and~\ref{thm:SWnonsmooth} use sieve methods, in particular the fundamental lemma of sieve theory (Theorem~\ref{thm:FLST}).

The main difficulty in proving Theorem~\ref{thm:main} is that the Golomb primes are sparse and that little is known about them beyond Erd\H{o}s' result, Theorem~\ref{thm:erdos}. Theorems~\ref{thm:SW} and~\ref{thm:SWnonsmooth} let us find subsets of shifted Golomb primes to which a proposition of the kind introduced by Croot~\cite{Croot2003} applies. Because the Golomb primes are sparse, the shifted Golomb primes in a single interval $(N^{0.997},N]$ have reciprocal sum  of order only $1/\log\log N$, so we have to combine about $\log\log N$ intervals. We also need a single modulus $P$, divisible by $b$, that divides the least common multiple of the set chosen in each interval. For this reason, our method does not extend Theorem~\ref{thm:main} to an arbitrary positive density subset of the Golomb primes.
\subsection*{Acknowledgements} The author thanks her advisors Giorgis Petridis and Paul Pollack for many helpful discussions and their constant encouragement. The author was supported by NSF DMS Award 2054214.

\subsection*{Use of AI statement}
The author used Claude Fable 5.1 and Claude Opus 5.5 to check the mathematical arguments for errors and gaps, and to proofread the paper. Apart from one suggestion by Fable 5.1, described below, all proof ideas are the author's own.

The author originally planned to use~\cite[Proposition~5.2]{LS24} in place of Proposition~\ref{prop:main}. While writing the paper, she found that the argument needs condition~(iii) of Proposition~\ref{prop:main} with a smaller value of $\eta$ than~\cite[Proposition~5.2]{LS24} allows. She asked Fable 5.1 whether the proof of~\cite[Proposition~5.2]{LS24} could be adapted to give the smaller value. It replied that the improvement does not follow directly, and suggested combining Bloom's treatment of the minor arcs~\cite{Bloom2025density} with the fact that the set $A$ consists of shifted primes. The author worked out the details; the resulting proof of Proposition~\ref{prop:main} is in Appendix~\ref{appendix}. The author has checked all proofs and takes full responsibility for the content of the paper.
\section{Notation and outline of the proof}\label{sec:notation}
\subsection*{Notation} We use the following notation.
\begin{itemize}
    \item We use $\mathcal{P}$ to denote the set of primes, $\mathcal{P}(x)$ to denote the set of all primes $\leq x$ and $\mathcal{P}(x;b,h)$ to denote the set of all primes $\leq x$ which lie in the congruence class $h\pmod{b}$.
    \item We use $\mathcal{G}$ to denote the set of Golomb primes,  $\mathcal{G}(x)$ to denote the set of all Golomb primes $\leq x$ and $\mathcal{G}(x;b,h)$ to denote the set of all Golomb primes $\leq x$ which lie in the congruence class $h\pmod{b}$.
    \item An $h$-shifted Golomb prime is a positive integer of the form $g-h$ with $g\in\mathcal{G}$. An integer $n$ is Golomb-prime-free if no Golomb prime divides $n$.
    \item We follow the convention that $(b,0)=b$ and that every number $n$ divides $0$.
    \item For a finite set $A$ of positive integers, let \[R(A):=\sum_{n\in A} \frac{1}{n}.\] For a prime power $q$, we write $q\pdiv n$ if $q\mid n$ and $(q,n/q)=1$, and we put \[A_{q}:=\{n\in A:\ q\pdiv n\}\qquad\text{and}\qquad R(A;q):=\sum_{n\in A_{q}} \frac{q}{n}.\] We write $\mathcal{Q}_{A}$ for the set of prime powers $q$ such that $q\pdiv n$ for some $n\in A$, and $[\mathcal{D}]$ for the least common multiple of a finite set $\mathcal{D}$ of integers. This notation follows Bloom~\cite{Bloom2025density} and is used mainly in Section~\ref{mainthm} and Appendix~\ref{appendix}.
    \item $\omega(n)$ is the number of distinct prime factors of $n$, and for a prime $l$, $\nu_{l}(n)$ is the exponent of $l$ in $n$, that is, the largest $k$ with $l^{k}\mid n$. We write $e(\theta)=e^{2\pi i\theta}$, and $\|\theta\|$ for the distance from $\theta$ to the nearest integer.
    \item $\log_{k}$ is the $k$-fold iterated logarithm, so $\log_{2}x=\log\log x$ and $\log_{3}x=\log\log\log x$, and $\operatorname{Li}(x):=\int_{2}^{x}\frac{dt}{\log t}$.
    \item Unless stated otherwise, $p$ and $r$ denote primes, $l$ denotes a prime that is not a Golomb prime, and $g$ and $q$ denote Golomb primes. In Section~\ref{mainthm} and Appendix~\ref{appendix}, $q$ denotes a prime power instead.
    \item We use standard asymptotic notation. For quantities $X$ and $Y$ with $Y>0$,
    we write $X=O(Y)$ or $X\ll Y$  if $|X|\leq CY$ for some constant $C>0$, and  $X\asymp Y$ or $X=\Theta(Y)$ if both $X\ll Y$ and $Y\ll X$. We write $X=o(Y)$ if $X/Y\to0$. Subscripts, as in $O_h$, $\ll_A$, show which parameters the implied constants depend on.
\end{itemize}
\subsection*{Outline of the proof of Theorem~\ref{thm:main}}

The proofs in~\cite{Bloom2025density,Bloom2025shifted,LS24} all rely on a proposition of the kind introduced by Croot~\cite{Croot2003}. Roughly speaking, it says the following. Let $A\subseteq(M,N]$ and $Q=\lcm(A)$. Then the equation
\[
  \frac{x}{Q}=\sum_{1\leq i\leq k}\frac{1}{n_i}
\]
has a solution in distinct $n_1,\dots,n_k\in A$, provided that the following conditions hold:

\begin{enumerate}
  \item Every $n\in A$ is power-smooth: every prime power
  dividing $n$ is at most $N^{0.99}$.
  \item Every $n\in A$ has ``about'' $\log\log n$ prime factors.
  \item Every prime power dividing $Q$ divides sufficiently
  many elements of $A$.
  \item  $x/Q$ is ``comparable'' to
  $\sum_{n\in A}1/n$.
\end{enumerate}
The precise version we use is Proposition~\ref{prop:main}.
Applying it to the $h$-shifted Golomb primes presents two
difficulties.

The first is to find subsets of shifted Golomb primes satisfying
the hypotheses of the proposition, as not much is known about
the arithmetic of Golomb primes beyond the estimates of
Erd\H{o}s~\cite{Erdos1961} recalled in
Theorem~\ref{thm:erdos}. We therefore prove two new results:
a Siegel--Walfisz-type theorem for Golomb primes in arithmetic
progressions (Theorem~\ref{thm:SW}) in Section~\ref{sec:SW}, and an upper bound for
the number of $h$-shifted Golomb primes that are not power-smooth (Theorem~\ref{thm:SWnonsmooth}) in Section~\ref{sec:SWnonsmooth}. Both proofs use the
sieve estimates recorded in Section~\ref{sec:SE}.

The second difficulty is that shifted Golomb primes are too sparse for a single application of the proposition to
suffice. In the earlier works~\cite{Croot2003,Bloom2025density,Bloom2025shifted,LS24},
it is enough to apply the proposition to a suitable set in
one interval $[M,N]$. Our Proposition~\ref{prop:main} requires
$A\subseteq(M,N]$ with $M\geq N^{0.997}$. However, the
reciprocal sum of the $h$-shifted Golomb primes in such an
interval is only of order $1/\log_{2}N$
(see Lemma~\ref{lem:RSA'}). For large $N$, this is too small
to represent a fixed positive fraction $a/b$.

Instead, we split the $h$-shifted Golomb primes into $J\asymp\log_{2}N$ consecutive intervals $(N_{j+1},N_{j}]$, where $N_{j+1}=N_{j}^{e^{-0.003}}$, and in the $j$th interval we find a set $A_{j}$ to which Proposition~\ref{prop:main} applies. Together, these sets have reciprocal sum greater than $1$. The key step is to choose an integer $P$ of size about $(\log_{2}N)^{2}$ such that $b\mid P$ and $P\mid\lcm(A_{j})$ for every $j$ (Lemma~\ref{lem:P}). We need this to apply Proposition~\ref{prop:main} to represent a suitable fraction $y_{j}/P$. Since $b\mid P$, we can choose the $y_{j}$ so that
\[
  \sum_{j=1}^{J}\frac{y_j}{P}=\frac{a}{b}.
\]
The intervals are disjoint, so adding these representations gives a representation of $a/b$ as a sum of distinct unit fractions with denominators of the form $g-h$.
 In Section~\ref{mainthm}, we make this argument precise and deduce
Theorem~\ref{thm:main} from Proposition~\ref{prop:main}. We prove Proposition~\ref{prop:main} in Appendix~\ref{appendix}.

\section{Some technical estimates}\label{sec:SE}
  This section collects some tools we use in the rest of the paper.
We first recall two elementary estimates of Mertens which we use several times:
\begin{equation}\label{mertens1}
  \sum_{p^k\leq X}\frac{1}{p^k}=\log\log X+c+O\Bigl(\frac{1}{\log X}\Bigr)
  \quad\text{for some constant } c,
\end{equation}
and
\begin{equation}\label{mertens2}
      \prod_{p\leq X}\Bigl(1-\frac{1}{p}\Bigr)^{-1}=\bigl(e^{\gamma}+o(1)\bigr)\log X.
\end{equation}
Here the sum in \eqref{mertens1} runs over prime powers $p^k$, and the product in \eqref{mertens2} runs over primes $p$.

  We also use the Bombieri--Vinogradov theorem~\cite{Bombieri1965, Vinogradov1965}. For $(a,q)=1$, let $\pi(x;q,a)$ be the number of primes $p\leq x$ with $p\equiv a\pmod q$, and let \[E(x;q,a):=\pi(x;q,a)-\frac{\operatorname{Li}(x)}{\phi(q)}.\] In the form given in~\cite[Chapter~28]{Davenport2000}, combined with partial summation, the theorem states that for every fixed $A>0$,
  \begin{equation}\label{eq:BV}
  \sum_{q\leq x^{1/2}/(\log x)^{A+8}}\ \max_{(a,q)=1}|E(x;q,a)|\ll_{A}\frac{x}{(\log x)^{A+3}}.
  \end{equation}

Next, we state a result for the sum of the reciprocals of the Golomb primes in an interval $(y,x]$. 
\begin{lemma}\label{lem:partialsum}
For $3\leq y<x$, we have
    \[\sum_{\substack{q\in\mathcal{G}\\ y< q\leq x}}\frac{1}{q}=(1+o(1))(\log_{3}x -\log_{3} y)+O\bigg(\frac{1}{\log y\log_{2} y}\bigg),\]
where the $o(1)$ term tends to $0$ as $y\to\infty$.
\end{lemma}
\begin{proof}
By partial summation and Theorem~\ref{thm:erdos}, we have
    \[
        \sum_{\substack{q\in\mathcal{G}\\ y< q\leq x}}\frac 1q = \frac{x(1+o(1))}{x\log x\log_{2} x} - \frac{y(1+o(1))}{y\log y\log_{2} y}+\int_{y}^{x}\frac{(1+o(1))}{u\log u\log_{2}u}\,du.
    \]
  Since $\displaystyle\int_{y}^{x}\frac{du}{u\log u\log_{2}u}=\log_{3}x-\log_{3}y$, this gives
    \[
         \sum_{\substack{q\in\mathcal{G}\\ y< q\leq x}}\frac 1q =(1+o(1))(\log_{3}x -\log_{3} y)+O\bigg(\frac{1}{\log y\log_{2} y}\bigg).\qedhere
    \]
\end{proof}

We also need three sieve estimates, which we use in Sections~\ref{sec:SW},~\ref{sec:SWnonsmooth} and in Appendix~\ref{appendix}.
\begin{theorem}[{Fundamental lemma of sieve theory;~\cite[Theorem~18.11]{Koukoulopoulos2019}}]\label{thm:FLST}
Let $\mathcal{A}$ be a finite set of positive integers, let $\mathcal{P}$ be a set of primes, and let $P(z)$ be the product of the primes $p\in\mathcal{P}$ with $p<z$. Let $S(\mathcal{A},\mathcal{P},z)=|\mathcal{A}\setminus\bigcup_{p\mid P(z)} \mathcal{A}_{p}|$, where $\mathcal{A}_{d}$ denotes the set of multiples of $d$ in $\mathcal{A}$. Suppose that the following conditions are satisfied:
\begin{enumerate}
    \item There exist a multiplicative function $w$, a parameter $X$ and a sequence of errors $r_{d}$ for $d\mid P(z)$ such that \[|\mathcal{A}_{d}|=\frac{{w(d)}}{d}X+r_{d}\quad\text{for all } d\mid P(z),\qquad\text{and}\qquad w(p)<p\quad\text{for all } p\mid P(z).\]
    \item There are constants $\kappa \geq 0 $ and $C\geq 0 $ such that
    \[
    \prod_{\substack{p\in \mathcal{P}\\y_{1}< p\leq y_{2}}}\left( 1-\frac{{w(p)}}{p} \right)^{-1}\leq\left( 1+\frac{C}{\log{y_{1}}} \right)\left( \frac{{\log y_{2}}}{\log y_{1}} \right)^{\kappa}
    \]
     uniformly for $\frac{3}{2}\leq y_{1}\leq y_{2}\leq z$. \end{enumerate}
     Then, uniformly for all $u\geq 1 $, we have \[S(\mathcal{A},\mathcal{P},z)=X\prod_{\substack{p<z\\p\in \mathcal{P}}} \left( 1-\frac{{w(p)}}{p} \right)(1+O_{\kappa,C}(u^{-u/2}))+O\left( \sum_{\substack{d\leq z^u\\d\mid P(z)}}|r_{d}|  \right).\]

\end{theorem}
We also use the following form of the fundamental lemma for the Selberg sieve, which is stated in \cite[Theorem~7.2]{HalberstamRichert1974}. This can be proved using Theorem \ref{thm:FLST}, but this form is more convenient.

\begin{theorem}[Selberg sieve]\label{thm:FLSS} Let $\mathcal{A}$ be a finite set of positive integers, let $\mathcal{P}$ be a set of primes, and let $P(z)$ be the product of the primes $p\in\mathcal{P}$ with $p<z$. Let $S(\mathcal{A},\mathcal{P},{z})=|\mathcal{A}\setminus\bigcup_{p\mid P(z)} \mathcal{A}_{p}|$. Suppose that $\displaystyle|\mathcal{A}_{d}|= \frac{w(d)}{d}X+R_{d}$ for all squarefree $d\mid P(z)$, where $w$ is a multiplicative function, and that the following conditions are satisfied:
\begin{enumerate}
    \item $\displaystyle 0\leq \frac{w(p)}{p}\leq 1-c$ for some absolute constant $c>0$;
    \item $|R_{d}|\leq w(d)$ for all squarefree $d\mid P(z)$;
    \item $\displaystyle\sum_{\substack{y\leq p\leq z\\p \in \mathcal{P}}} \frac{{w(p)\log p}}{p}\leq \kappa \log\left( \frac{z}{y} \right)+A_{2}$ for all $2\leq y\leq z$, where $\kappa$ and $A_{2}$ are constants.
\end{enumerate}
If, in addition, $X\geq z$, then
\[S(\mathcal{A},\mathcal{P},z)=X\prod_{\substack{p<z\\p\in \mathcal{P}}}\left( 1-\frac{w(p)}{p} \right)  \left( 1+O\left( \exp\left( -\frac{{\log X}}{2\log z} \right) \right) \right).\]

\end{theorem}
In Theorems~\ref{thm:FLST} and~\ref{thm:FLSS}, $\mathcal{A}_{p}$ is the set of $n\in\mathcal{A}$ with $p\mid n$, so we remove the multiples of each prime $p\mid P(z)$. In our applications we instead remove the elements of $\mathcal{A}$ that lie in some given residue classes modulo $p$. The following observation shows that both theorems still apply in this setting. It extends the argument in the proof of~\cite[Lemma~2.4]{Stucky2025}.
\begin{observation}\label{lem:residuesieve}
Let $\mathcal{A}$ be a finite set of positive integers, and let $\mathcal{P}$ and $P(z)$ be as in Theorem~\ref{thm:FLST}. For each prime $p\mid P(z)$, suppose that we remove from $\mathcal{A}$ the elements that lie in $w(p)\geq 1$ given residue classes modulo $p$. For $d\mid P(z)$, let $\mathcal{B}_{d}$ be the set of $n\in\mathcal{A}$ that lie in a removed residue class modulo every prime $p\mid d$. Let $k=\max_{p\mid P(z)}w(p)$, and write the classes removed modulo $p$ as $a_{p,1},a_{p,2},\ldots,a_{p,k}$, repeating a class if $w(p)<k$. By the Chinese remainder theorem, we can choose integers $a_{1},\ldots,a_{k}\leq 0$ such that $a_{i}\equiv a_{p,i}\pmod{p}$ for every $i$ and every prime $p\mid P(z)$. Define
\[
F(T)=\prod_{i=1}^{k}(T-a_{i})\qquad\text{and}\qquad \tilde{\mathcal{A}}=\{F(n):\ n\in\mathcal{A}\}.
\]
Since each factor $n-a_{i}$ is positive and increases with $n$, $F$ is one-to-one on $\mathcal{A}$. For a prime $p\mid P(z)$, we have $p\mid F(n)$ if and only if $n\equiv a_{p,i}\pmod{p}$ for some $i$, that is, if and only if $n$ lies in a removed class modulo $p$. Since every $d\mid P(z)$ is squarefree, it follows that $|\mathcal{B}_{d}|$ equals the number of elements of $\tilde{\mathcal{A}}$ divisible by $d$. In particular, $|\mathcal{A}\setminus\bigcup_{p\mid P(z)}\mathcal{B}_{p}|$ equals the number of elements of $\tilde{\mathcal{A}}$ coprime to $P(z)$. Applying Theorems~\ref{thm:FLST} and~\ref{thm:FLSS} to $\tilde{\mathcal{A}}$ therefore shows that they still hold with the sets $\mathcal{B}_{d}$ in place of the multiples of $d$ in $\mathcal{A}$.
\end{observation}
In the rest of the paper we use Theorems \ref{thm:FLST} and \ref{thm:FLSS} in this form. We abuse notation and write $\mathcal{A}_{d}$ and $\mathcal{S}(\mathcal{A},\mathcal{P},z)$ instead of $\mathcal{B}_{d}$ and $\mathcal{S}(\mathcal{\tilde{A}},\mathcal{P},z)$ respectively.
The following lemma follows from a result in Diamond--Halberstam~\cite[Theorem~5.6]{DiamondHalberstam2008}.
\begin{lemma}\label{lem:sievebound} 
     Let $x$ be large and let $z=\frac{x^{1/4}}{(\log{x})^{100}}$. For each prime $r\leq z$, remove from the primes $p\leq x$ those lying in $w(r)$ given residue classes modulo $r$, each coprime to $r$, where $w(r)\leq \min\{2,r-2\}$. Then, uniformly in the choice of these residue classes, the number of remaining primes is at most
\[
\frac{2e^{2\gamma}x(1+o(1))}{\log{x}}\prod_{r\leq z}\left(1-\frac{w(r)}{r-1}\right).
\]
 \end{lemma}
      \begin{proof}
     Let $\mathcal{A}$ be the set of primes $p\leq x$, let $\mathcal{P}$ be the set of all primes and let $z=\frac{x^{1/4}}{(\log{x})^{100}}$. Let $\mathcal{A}_{p}$ be the set of elements of $\mathcal{A}$ that lie in $w(p)$ given residue classes modulo $p$. By the Chinese remainder theorem, for squarefree $n\mid P(z)$, $\mathcal{A}_{n}$ is the set of primes $\leq x$ that lie in $w(n)$ given residue classes modulo $n$, all coprime to $n$, so
     \[
     |\mathcal{A}_{n}|=\frac{w(n)}{\phi(n)}\operatorname{Li}(x)+r_{\mathcal{A}}(n),
     \]
     where $r_{\mathcal{A}}(n)$ is the sum of $E(x;n,a)$ over these $w(n)$ residue classes $a$, with $E(x;n,a)$ as in equation \eqref{eq:BV}. Since $w(n)\leq 2^{\omega(n)}$, it follows that
     \[
     |r_{\mathcal{A}}(n)|\leq 2^{\omega(n)}E(x;n),\qquad\text{where}\quad E(x;n):=\max_{(a,n)=1}|E(x;n,a)|.
     \]
     We check the conditions and apply~\cite[Theorem~5.6]{DiamondHalberstam2008} to the set $\mathcal{A}$, with $X=\operatorname{Li}(x)$, $q=1$, $u=1$ and $z=x^{1/4}/(\log x)^{100}$. This gives
     \[
     S(\mathcal{A},\mathcal{P},z)\leq \operatorname{Li}(x)\prod_{r\leq z}\left(1-\frac{w(r)}{r-1}\right)\bigg\{2e^{2\gamma}+O\left( \frac{\log_{2}x}{\log x} \right)\bigg\}+\sum_{\substack{n<z^2\\n\mid P(z)}}3^{\omega(n)}|r_{\mathcal{A}}(n)|.
     \]

     We bound the last sum using the Cauchy--Schwarz inequality:
     \[
     \sum_{\substack{n \mid P(z)\\n<z^2}}3^{\omega(n)}|r_{\mathcal{A}}(n)|\leq \sum_{\substack{n \mid P(z)\\n<z^2}}6^{\omega(n)}E(x;n)\leq \bigg(\sum_{\substack{n \mid P(z)\\n<z^2}}36^{\omega(n)}E(x;n)\bigg)^{1/2}\bigg(\sum_{\substack{n \mid P(z)\\n<z^2}}E(x;n)\bigg)^{1/2}.
     \]
     Since $z^{2}=x^{1/2}/(\log x)^{200}$, equation \eqref{eq:BV} gives
     \[
     \sum_{\substack{n \mid P(z)\\n<z^2}}E(x;n)\leq\sum_{n\leq x^{1/2}/(\log x)^{200}}\ \max_{(a,n)=1}|E(x;n,a)|\ll\frac{x}{(\log x)^{195}},
     \]
     so the second factor is $O\big(x^{1/2}/(\log x)^{97}\big)$. For the first factor, we use the trivial bound $E(x;n)\ll x/\phi(n)$, which gives
     \[
    \bigg( \sum_{\substack{n \mid P(z)\\n<z^2}}36^{\omega(n)}E(x;n)\bigg)^{1/2}\ll x^{1/2}\bigg(\sum_{n\mid P(z)}\frac{36^{\omega(n)}}{\phi(n)}\bigg)^{1/2}\ll x^{1/2}\bigg(\prod_{p\leq z}\left(1+\frac{36}{p-1}\right)\bigg)^{1/2}\ll x^{1/2}(\log x)^{18}.
     \]
     So the total error is $O(x/(\log x)^{79})$, which is negligible, since $\prod_{r\leq z}\big(1-\frac{w(r)}{r-1}\big)\gg(\log x)^{-2}$ and $\operatorname{Li}(x)=(1+o(1)) x/\log x$.
    Therefore, we obtain
     \[
      S(\mathcal{A},\mathcal{P},z)\leq\frac{2e^{2\gamma}x(1+o(1))}{\log x} \prod_{r\leq z}\left(1-\frac{w(r)}{r-1}\right),
     \]
      which is the claimed bound.
 \end{proof}
\section{Golomb primes in arithmetic progressions: proof of Theorem~\ref{thm:SW}} \label{sec:SW}
In this section we prove Theorem~\ref{thm:SW}, which we restate here.
\mtheoremtwo*

 Recall that $\mathcal{P}$ is the set of all primes, $\mathcal{P}(x;b,h)$ is the set of primes $p\leq x$ with $p\equiv h\pmod{b}$ and $\mathcal{G}(x;b,h)$ is the set of Golomb primes $g\leq x$ with $g \equiv h\pmod{b}$. Unless stated otherwise, $r$ denotes a prime, $l$ a prime that is not a Golomb prime, and $g$ and $q$ Golomb primes. 
 
To estimate $|\mathcal{G}(x;b,h)|$, we start with $\mathcal{P}(x;b,h)$ and remove the primes $p\equiv 1\pmod q$ for every Golomb prime $q$. For the Golomb primes $q\leq z_{0}:=\exp\big(\frac{\log x}{3\log_{3}x}\big)$, we do this with the fundamental lemma of sieve theory (Theorem~\ref{thm:FLST}). We then show that the remaining Golomb primes $z_{0}<q\leq x$ remove a number of primes that is small compared with the main term.

    Let $x$ be large and let $1\leq b\leq (\log x)^A$. We apply Theorem~\ref{thm:FLST} with
    \[
    \mathcal{A}=\mathcal{P}(x;b,h),\qquad X=\frac{\operatorname{Li}(x)}{\phi(b)}, \qquad \mathcal{P}=\{g\in \mathcal{G}:g\nmid b\}.
    \]
       Further, let $u=\log_{3}x$ and $z=z_{0}=\exp\left( \frac{{\log x}}{3\log_{3}x} \right)$, and let $E(x;q,a)$ be as in Section~\ref{sec:SE}. For squarefree $d\mid P(z_{0})$ we have $(d,b)=1$, so by the Chinese remainder theorem there is a unique residue class $h_{d}\pmod{db}$ with $h_{d}\equiv h\pmod{b}$ and $h_{d}\equiv 1 \pmod{d}$. Since $(h,b)=1$, this class is coprime to $db$, and \[|\mathcal{A}_{d}|=\frac{\operatorname{Li}(x)}{\phi(d)\phi(b)}+E(x;db,h_{d}).\]
    Thus, with $\displaystyle \frac{w(g)}{g}=\frac{1}{g-1}$ and $r_{d}=E(x;db,h_{d})$, the conditions of Theorem~\ref{thm:FLST} are satisfied with $\kappa=1$, and as $u^{-u/2}\leq(\log_{2}x)^{-2}$ we obtain
    \begin{equation}\label{eq:SAPz0}
    S(\mathcal{A},\mathcal{P},z_{0})=\frac{\operatorname{Li}(x)}{\phi(b)}\prod_{\substack{g\in\mathcal{G}(z_{0})\\g\nmid b}}\left( 1-\frac{1}{g-1} \right)\Big(1+O\big((\log_{2}x)^{-2}\big)\Big)+O\left( \sum_{\substack{d\leq z_{0}^u\\d\mid P(z_{0})}}|E(x;db,h_{d})|  \right).
   \end{equation}
      We now bound the error term. As $d\leq z_{0}^{u}=x^{1/3}$, we have $db\leq x^{1/3}(\log x)^{A}\leq x^{1/2}/(\log x)^{A+8}$, so equation \eqref{eq:BV} gives
      \[
       S(\mathcal{A},\mathcal{P},z_{0})=(1+o(1))\frac{\operatorname{Li}(x)}{\phi(b)}\prod_{\substack{g\in\mathcal{G}(z_{0})\\g\nmid b}}\left( 1-\frac{1}{g-1} \right)+O\left( \frac{x}{(\log x)^{A+3}} \right).
      \]

   Since $\phi(b)\leq(\log x)^A$, the error term is negligible in comparison to the main term. By Theorem~\ref{thm:erdos} and  $z_{0}=\exp\left( \frac{{\log x}}{3\log_{3}x} \right)$, we have
    \[
    \prod_{\substack{g \in \mathcal{G}(z_{0})\\g\nmid b}}\left( 1-\frac{1}{g-1} \right)=\frac{1+o(1)}{\log_{2}x}\prod_{\substack{g\in \mathcal{G}\\g\mid b}}\left( 1-\frac{1}{g-1} \right)^{-1}.
    \]
   Since $\operatorname{Li}(x)=(1+o(1))\frac{x}{\log x}$, we obtain \[S(\mathcal{A},\mathcal{P},z_{0})=(1+o(1))\frac{x}{\phi(b)\log x\log_{2}x}\prod_{\substack{g\in \mathcal{G}\\g\mid b}}\left( 1-\frac{1}{g-1} \right)^{-1}.\]

    Note that $S(\mathcal{A},\mathcal{P},z_{0})$ counts the primes $p\leq x$ with $p\equiv h\pmod{b}$ such that $p-1$ is not divisible by any Golomb prime $g\leq z_{0}$ with $g\nmid b$. Hence every element of $\mathcal{G}(x;b,h)$ is counted by $S(\mathcal{A},\mathcal{P},z_{0})$.
     On the other hand, to bound $|\mathcal{G}(x;b,h)|$ from below, we bound the number of primes counted by $S(\mathcal{A},\mathcal{P},z_{0})$ that are not Golomb primes. Define
      \[B(x;b,h; g):=\#\{p\leq x:\ p\equiv h\pmod{b},\ p\equiv 1\pmod{g},\ p\not\equiv 1\pmod{q}\text{ for all } q\in \mathcal{G}(z_{0})\}.\]

    Every prime $p\neq 2$ that is counted by $S(\mathcal{A},\mathcal{P},z_{0})$ but is not a Golomb prime is counted by $B(x;b,h;g)$ for at least one Golomb prime $g\in(z_{0},x]$. Therefore
   \begin{equation}\label{eq:Gxhb}
       S(\mathcal{A},\mathcal{P},z_{0})-1-\sum_{\substack{g\in\mathcal{G}\\ z_{0}<g\leq x}}B(x;b,h; g)\leq|\mathcal{G}(x;b,h)|\leq S(\mathcal{A},\mathcal{P},z_{0}).
   \end{equation}
     We split the range of $g$ into two subintervals at $\displaystyle t=x\exp\left(-\frac{\log{x}}{(\log_{2}{x})^{4/5}}\right)$ and bound the sum over $g$ in equation~\eqref{eq:Gxhb} using the following two lemmas.
   \begin{lemma}\label{lem:errorzt}
For $\displaystyle z_{0}=\exp\left( \frac{{\log x}}{3\log_{3}x} \right)$ and $\displaystyle t=x\exp\left(-\frac{\log{x}}{(\log_{2}{x})^{4/5}}\right),$ we have
\[\sum_{\substack{g\in\mathcal{G}\\ z_{0}< g\leq t}}B(x;b,h;g)\ll\frac{x\log_4x}{\phi(b)\log x\,(\log_2x)^{6/5}}\prod_{\substack{q\mid b\\q\in\mathcal{G}}}\Big(1-\frac{1}{q-1}\Big)^{-1}.\]
\end{lemma}
\begin{lemma}\label{lem:errortx}
For $\displaystyle t=x\exp\left(-\frac{\log{x}}{(\log_{2}{x})^{4/5}}\right),$ we have
\[\sum_{\substack{g\in\mathcal{G}\\ t<g\leq x}}B(x;b,h;g)\ll\frac{x(\log_3x)^3}{\phi(b)\log x\,(\log_2x)^{2}}.\]
\end{lemma}
We first assume the lemmas and deduce Theorem~\ref{thm:SW} from them.
By Lemmas~\ref{lem:errorzt} and~\ref{lem:errortx}, we have
\[
    \sum_{\substack{g\in\mathcal{G}\\ z_{0}<g\leq x}}B(x;b,h;g)\ll_{A}\frac{x\log_4x}{\phi(b)\log x\,(\log_2x)^{6/5}}\prod_{\substack{q\mid b\\q\in\mathcal{G}}}\Big(1-\frac{1}{q-1}\Big)^{-1}+\frac{x(\log_3x)^3}{\phi(b)\log x\,(\log_2x)^{2}}.
\]
Note that $(\log_3x)^3/(\log_2x)^{2}=o\big(\log_4x/(\log_2x)^{6/5}\big)$, so the first term dominates and
\[
    \sum_{\substack{g\in\mathcal{G}\\ z_{0}<g\leq x}}B(x;b,h;g)\ll_{A}\frac{x\log_4x}{\phi(b)\log x\,(\log_2x)^{6/5}}\prod_{\substack{q\mid b\\q\in\mathcal{G}}}\Big(1-\frac{1}{q-1}\Big)^{-1}.
\]
This is smaller than the main term. Combining this with \eqref{eq:Gxhb} and the asymptotic formula for $S(\mathcal{A},\mathcal{P},z_{0})$, we obtain
\[|\mathcal{G}(x;b,h)|=(1+o_{A}(1))\frac{x}{\phi(b)\log x\log_{2}x}\prod_{\substack{g\mid b\\g\in \mathcal{G}}}\left( 1-\frac{1}{g-1} \right)^{-1}. \]
This proves Theorem~\ref{thm:SW}, assuming Lemmas~\ref{lem:errorzt} and~\ref{lem:errortx}, which we now prove.

    \begin{proof}[Proof of Lemma~\ref{lem:errorzt}]
        Recall that \[B(x;b,h; g)=\#\{p\leq x:\ p\equiv h\pmod{b},\ p\equiv 1 \pmod{g},\ p\not\equiv 1 \pmod{q}\text{ for all } q\in \mathcal{G}(z_{0})\}.\]
        Fix a Golomb prime $g\in(z_{0},t]$. Since $g>z_{0}>b$, we have $(g,b)=1$. By the Chinese remainder theorem, there is a unique integer $\tilde{h}\in[0,gb)$ with $\tilde{h}\equiv h\pmod{b}$ and $\tilde{h}\equiv 1 \pmod{g}$. Then every prime $p$ counted by $B(x;b,h;g)$ can be written as
        \[
        p= gbm+\tilde{h}, \qquad\text{with }0\leq m\leq \frac{x}{gb}.
        \]
        The value $m=0$ occurs only for $p=\tilde{h}$, so these primes contribute at most $|\mathcal{G}(t)|\leq t$ to the sum in Lemma~\ref{lem:errorzt}. This is negligible compared to the bound we are aiming to prove; hence it suffices to count the primes $p=gbm+\tilde{h}$ with $1\leq m\leq x/(gb)$.
        Note that $(\tilde{h},gb)=1$ and that no Golomb prime $q\leq z_{0}$ divides $(\tilde{h}-1,gb)$. If $p$ is counted by $B(x;b,h;g)$, then $p\not\equiv 0 \pmod{r}$ for all $r\in\mathcal{P}$ with $r\not=p$, and $p\not\equiv 1 \pmod{q}$ for all $q\in\mathcal{G}(z_{0})$. We now translate these congruence conditions on $p$ into congruence conditions on $m$.
        \begin{itemize}
            \item For all $r\in\mathcal{P}$, if $r\nmid gb$, then $p= gbm+\tilde{h}\not\equiv 0 \pmod{r}$ implies $\displaystyle m\not\equiv-\frac{\tilde{h}}{gb}\pmod{r}$.
            \item For all $q\in\mathcal{G}(z_{0})$, if $q\nmid gb$, then $p= gbm+\tilde{h}\not\equiv 1 \pmod{q}$ implies $\displaystyle m\not\equiv-\frac{\tilde{h}-1}{gb}\pmod{q}$.
        \end{itemize}
        Note also that $\displaystyle-\frac{\tilde{h}}{gb}$ and $\displaystyle-\frac{\tilde{h}-1}{gb}$ are distinct residue classes modulo every prime not dividing $gb$.
        We apply Theorem~\ref{thm:FLSS} to the set $\mathcal{A}$ of positive integers $m\leq x/(gb)$, removing the residue classes above, with $\kappa=2$, $\displaystyle X=\frac{x}{gb}$ and $ \displaystyle z_{1}=\exp\left( \frac{{\log x}}{2(\log_{2}x)^{4/5}} \right)$. For a prime $l$ that is not a Golomb prime, \[|\mathcal{A}_{l}|=\frac{x}{gbl}+O(1)\quad\text{ if } l\nmid b, \text{ and } |\mathcal{A}_{l}|=0 \text{ otherwise.}\]
        For a Golomb prime $q$,\[|\mathcal{A}_{q}|=\frac{{2x}}{gbq}+O(1)\quad\text{ if } q\nmid b, \text{ and } |\mathcal{A}_{q}|=0 \text{ otherwise.}\]
        Thus $w(l)=1$ for $l\in\mathcal{P}\setminus\mathcal{G}$ with $l\nmid b$, $w(q)=2$ for $q\in\mathcal{G}$ with $q\nmid b$, and $w(r)=0$ for all other $r\in\mathcal{P}$. We extend $w$ multiplicatively to all squarefree $d\mid P(z_{1})$.
        Since $g\leq t$, we have $X=x/(gb)\geq x/(tb)\geq z_{1}$, so Theorem~\ref{thm:FLSS} applies, and $1+O\big(\exp(-\frac{\log X}{2\log z_{1}})\big) =O(1)$. Therefore
        \[
            S(\mathcal{A},\mathcal{P},z_{1})\ll\frac{x}{gb}\prod_{\substack{l \nmid b\\l\in \mathcal{P}\setminus \mathcal{G}\\ l \leq z_{1}}}\left( 1-\frac{1}{l} \right) \prod_{\substack{q \nmid b\\q\in \mathcal{G}\\ q \leq z_{1}}}\left( 1-\frac{2}{q} \right).
        \]
        Since $1-\frac{2}{q}=\left(1-\frac{1}{q}\right)\left(1-\frac{1}{q-1}\right)$ and every prime factor of $b$ is less than $z_{1}$, it follows that
        \[
            S(\mathcal{A},\mathcal{P},z_{1})\ll\frac{x}{gb} \prod_{\substack{r\in \mathcal{P}\\ r \leq z_{1}}}\left( 1-\frac{1}{r} \right)   \prod_{\substack{q\in \mathcal{G}\\ q \leq z_{1}}}\left( 1-\frac{1}{q-1} \right)  \prod_{\substack{r\mid b\\r\in \mathcal{P}}}\left( 1-\frac{1}{r} \right)^{-1}  \prod_{\substack{q\mid b\\q\in \mathcal{G}}}\left( 1-\frac{1}{q-1} \right)^{-1}.
        \]
            By Mertens' theorem \eqref{mertens2}, Theorem~\ref{thm:erdos} and the definition of Euler's totient function $\phi$, we get
            \[
             S(\mathcal{A},\mathcal{P},z_{1})\ll \frac{x}{g\phi(b)\log{z_{1}} \log_{2}{z_{1}}}\prod_{\substack{q\mid b\\q\in \mathcal{G}}}\left( 1-\frac{1}{q-1} \right)^{-1}.
            \]
            Since $\log z_{1}=\frac{\log x}{2(\log_{2}x)^{4/5}}$ and $\log_{2}z_{1}=\log_{2}x+O(\log_{3}x)$, we get
            \[
             S(\mathcal{A},\mathcal{P},z_{1})\ll \frac{x}{g\phi(b)\log{x} (\log_{2}x)^{1/5}}\prod_{\substack{q\mid b\\q\in \mathcal{G}}}\left( 1-\frac{1}{q-1} \right)^{-1}.
            \]
        As $B(x;b,h;g)\leq  S(\mathcal{A},\mathcal{P},z_{1})$ for each such $g$, this bound gives
        \[
            \sum_{\substack{g\in \mathcal{G}\\z_{0}<g\leq t}}B(x;b,h;g) \ll \frac{x}{\phi(b)\log{x} (\log_{2}x)^{1/5}}\prod_{\substack{q\mid b\\q\in \mathcal{G}}}\left( 1-\frac{1}{q-1} \right)^{-1}\sum_{\substack{g\in \mathcal{G}\\z_{0}<g\leq t}} \frac{1}{g}.
        \]
        By Lemma~\ref{lem:partialsum}, it follows that
        \begin{align} \label{eq:Bzt}
            \sum_{\substack{g\in \mathcal{G}\\z_{0}<g\leq t}}B(x;b,h;g) &\ll \frac{x}{\phi(b)\log{x} (\log_{2}x)^{1/5}}\prod_{\substack{q\mid b\\q\in \mathcal{G}}}\left( 1-\frac{1}{q-1} \right)^{-1} \left( \log_{3}t-\log_{3}z_{0}+O\left(\frac{1}{\log z_{0}\log_{2}z_{0}}\right) \right).
        \end{align}
        By the definitions of $t$ and $z_{0}$, \[\log_{3}t=\log_{3}x+O\left( \frac{1}{(\log_{2} x)^{9/5}} \right),\quad \log_{3}z_{0}=\log_{3}x+O\left( \frac{{\log_{4}x}}{\log_{2}x} \right), \quad \frac{1}{\log z_{0}\log_{2}z_{0}}=O\left(\frac{{\log_{3}x\log_{2}x}}{\log x}\right).\]
       So the factor in parentheses in \eqref{eq:Bzt} is $\ll\log_{4}x/\log_{2}x$, and inserting this into \eqref{eq:Bzt}, we obtain
        \[
            \sum_{\substack{g\in \mathcal{G}\\z_{0}<g\leq t}}B(x;b,h;g) \ll \frac{x\log_{4}x}{\phi(b)\log{x} (\log_{2}x)^{6/5}}\prod_{\substack{q\mid b\\q\in \mathcal{G}}}\left( 1-\frac{1}{q-1} \right)^{-1}.
        \]
        This completes the proof of Lemma~\ref{lem:errorzt}.
        \end{proof}
        \begin{proof}[Proof of Lemma~\ref{lem:errortx}] Let $g\in(t,x]$ be a Golomb prime, and let $p$ be a prime counted by $B(x;b,h;g)$. Then \[p=gm+1\] for some positive integer $m$, and $m<x/t$ since $p\leq x$ and $g>t$. Thus $m$ is small, so instead of counting the primes $p$ for each $g$, we fix $m$ and count the Golomb primes $g\in(t,x]$ for which $gm+1$ is a prime counted by $B(x;b,h;g)$.
        Let $\mathcal{M}$ be the set of positive integers $m$ for which there is at least one such Golomb prime $g$. We first list the conditions on $m$.
        \begin{itemize}
            \item Since $p=gm+1$ is an odd prime and $g$ is odd, $m$ is even.
           \item If $q\in\mathcal{G}(z_{0})$, then $p=gm+1\not\equiv 1\pmod{q}$ gives $q\nmid m$. Every prime factor of $m$ is less than $x/t<z_{0}$, so $m$ has no Golomb prime factors.
          
            \item $p=gm+1\equiv h\pmod{b}$ gives $gm\equiv h-1\pmod{b}$. Since $(g,b)=1$, this forces \[(m,b)=(h-1,b)=:D.\] Writing $m=m'D$, $h-1=h'D$ and $b=b'D$ and dividing by $D$, we get \[gm'\equiv h'\pmod{b'}, \text{ which implies } g\equiv \frac{h'}{m'}\pmod{b'}.\]
        \end{itemize}
        Thus
        \begin{equation}\label{eq:defM}
             \mathcal{M}\subseteq\left\{ 1\leq m\leq \frac{x}{t}:\ 2\mid m,\ q\nmid m \text{ for all }q\in \mathcal{G},\ (m,b)=D \right\}.
        \end{equation}

        Let $z_{2}=x^{1/2}$. Next, we list the conditions on $g$ for a fixed $m \in \mathcal{M}$.
        \begin{itemize}
            \item As $g$ is a prime, we have $g\not\equiv 0 \pmod{r}$ for all primes $r\leq z_{2}$.
            \item As $g$ is a Golomb prime, we have $g\not\equiv 1 \pmod{q}$ for all Golomb primes $q\leq z_{2}$.
            \item As $p=gm+1$ is a prime and $p>z_{2}$, we have $p\not\equiv 0 \pmod{r}$ for all primes $r\leq z_{2}$. If $r\nmid m$, this implies $\displaystyle g\not\equiv-\frac{1}{m}\pmod{r}$.
        \end{itemize}
         The conditions modulo $3$ give one more restriction on $m$. Suppose that $3\nmid b'(m+1)$. Since $3\in \mathcal{G}$, we have $3\nmid m$, and so $m\equiv 1 \pmod{3}.$ Then $g$ would have to avoid the three classes  $0,$ $1$ and $-1/m\equiv 2 \pmod{3},$ which is impossible. Hence every $m\in \mathcal{M}$ satisfies $3\mid b'(m+1)$.
         
        For each $m\in \mathcal{M}$, let \(\mathcal{A}_{m}\) be the set of positive integers $g\leq x/m$ with \[g\equiv\frac{h'}{m'}\pmod{b'} \text{  so that  } |\mathcal{A}_{m}|=\frac{x}{mb'}+O(1).\]
        We now apply Theorem~\ref{thm:FLSS} to $\mathcal{A}=\mathcal{A}_{m}$, removing the residue classes above, with $\kappa=3$ and $z=z_{2}$. By the conditions on $g$ above, we can take
        \begin{itemize}
            \item $w(l)=2$ for every $l\in \mathcal{P}\setminus\mathcal{G}$ with $l\nmid m$ and $l \nmid b$
            \item $w(l)=1$ for every $l\in \mathcal{P}\setminus\mathcal{G}$ with $l\mid m$ and $l \nmid b'$
            \item $w(q)=3$ for every $q\in \mathcal{G}$ with $q \nmid b'$ and $q\nmid (m+1)$
            \item $w(q)=2$ for every $q\in \mathcal{G}$ with $q\nmid b$ and $q\mid m+1,$ since then $-1/m\equiv 1 \pmod{q}$, so two of three classes coincide,
        \end{itemize}
        and $w(r)=0$ for every prime $r \mid b'$.
         
        We extend $w$ multiplicatively to all squarefree $d\mid P(z_{2})$. Since $x/mb'>x^{0.98}$, Theorem \ref{thm:FLSS} gives
        \begin{align*}
            S(\mathcal{A}_{m},\mathcal{P},z_{2})&\ll\frac{x}{mb'}\prod_{\substack{l\nmid mb\\l\in \mathcal{P}\setminus \mathcal{G}\\l \leq x^{1/2}}}\left( 1-\frac{2}{l} \right) \prod_{\substack{l\mid m\\l\nmid b'\\l\in \mathcal{P}\setminus \mathcal{G}\\l \leq x^{1/2}}}\left( 1-\frac{1}{l} \right) \prod_{\substack{q \nmid b(m+1)\\q\in  \mathcal{G}\\q\leq x^{1/2}}}\left( 1-\frac{3}{q} \right)\prod_{\substack{q \nmid b\\q\mid(m+1)\\q\in  \mathcal{G}\\q\leq x^{1/2}}}\left( 1-\frac{2}{q} \right).
        \end{align*} Writing $1-\frac{2}{l}=\left( 1-\frac{1}{l} \right)\left( 1-\frac{1}{l-1} \right)$ for $l>2$ and $1-\frac{3}{q}=\left( 1-\frac{2}{q} \right)\left( 1-\frac{1}{q-2} \right)$ for $q>3$, we get
        \[S(\mathcal{A}_{m},\mathcal{P},z_{2})  \ll\frac{x}{mb'}\prod_{\substack{l\nmid b'\\ l\in\mathcal{P}\setminus\mathcal{G}\\ l\le x^{1/2}}}\Big(1-\frac1l\Big)\prod_{\substack{l\nmid m\\ l\nmid b'\\ l\in\mathcal{P}\setminus\mathcal{G}\\ 2<l\le x^{1/2}}}\Big(1-\frac1{l-1}\Big)\prod_{\substack{q\nmid b'\\ q\in\mathcal{G}\\ q\le x^{1/2}}}\Big(1-\frac2q\Big)\prod_{\substack{q\nmid b'\\ q\nmid m+1\\ q\in\mathcal{G}\\ 3< q\le x^{1/2}}}\Big(1-\frac1{q-2}\Big)\]
        Collecting the terms, we get
        \begin{align*}
           S(\mathcal{A}_{m},\mathcal{P},z_{2})  &\ll\frac{x}{mb'}\prod_{\substack{r\in\mathcal{P}\\ 2<r\le x^{1/2}}}\Big(1-\frac2r\Big)\prod_{\substack{r\mid b'\\r\in\mathcal{P}\\ r>2}}\Big(1-\frac2r\Big)^{-1}\prod_{\substack{l\mid m\\ l\nmid b'\\ l\in\mathcal{P}\setminus\mathcal{G}\\ 2<l\le x^{1/2}}}\Big(1-\frac1{l-1}\Big)^{-1}\\
            &\quad\prod_{\substack{q\in\mathcal{G}\\ 3<q\le x^{1/2}}}\Big(1-\frac1{q-2}\Big)\prod_{\substack{q\mid b'\\q\in\mathcal{G}\\ q>3}}\Big(1-\frac1{q-2}\Big)^{-1}\prod_{\substack{ q\nmid b'\\ q\mid m+1\\q\in \mathcal{G}\\ q>3}}\Big(1-\frac1{q-2}\Big)^{-1}.
            \end{align*}
            Every integer $1<n\leq x$ has $\ll \log x$ prime factors, so  Theorem \ref{thm:erdos} implies \[ \prod_{\substack{q\mid n,\,q >3\\q\in \mathcal{G}}} \left( 1-\frac{1}{q-2} \right)^{-1}\ll \log_{3}(x).
            \]
            Using this with $n=b'$ and $n=m+1$ for the last two products, we get
           
             \begin{align*}
           S(\mathcal{A}_{m},\mathcal{P},z_{2})  &\ll_{A}\frac{x(\log_{3}x_{})^2}{mb'}\prod_{\substack{r\in\mathcal{P}\\ 2<r\le x^{1/2}}}\Big(1-\frac2r\Big)\prod_{\substack{r\mid b'\\r\in\mathcal{P}\\ r>2}}\Big(1-\frac2r\Big)^{-1}\prod_{\substack{l\mid m\\ l\nmid b'\\ l\in\mathcal{P}\setminus\mathcal{G}\\ 2<l\le x^{1/2}}}\Big(1-\frac1{l-1}\Big)^{-1}\\
            &\quad\prod_{\substack{q\in\mathcal{G}\\ 3<q\le x^{1/2}}}\Big(1-\frac1{q-2}\Big).
            \end{align*}
        Consequently, \begin{align}\label{eq:errortx}
        \sum_{\substack{g\in \mathcal{G}\\t<g\leq x}} B(x;b,h;g)\leq\sum_{\substack{m\in \mathcal{M}\\1\leq m\leq x/t}}S(\mathcal{A}_{m},\mathcal{P},z_{2}) &\ll_A\frac{x(\log_{3}x)^2}{b'} \prod_{\substack{r\in\mathcal{P}\\ 2<r\le x^{1/2}}}\Big(1-\frac2r\Big)\prod_{\substack{r\mid b'\\r\in\mathcal{P}\\ r>2}}\Big(1-\frac2r\Big)^{-1}\prod_{\substack{q\in\mathcal{G}\\ 3<q\le x^{1/2}}}\Big(1-\frac1{q-2}\Big) \notag\\
        &\quad\sum_{\substack{m\in \mathcal{M}\\1\leq m\leq x/t}}
        \frac{1}{m} \prod_{\substack{l\mid m\\l \nmid b'\\l\in \mathcal{P}\setminus \mathcal{G}\\2<l \leq x^{1/2}}}\left( 1-\frac{1}{l-1} \right)^{-1}.\end{align}
        Let $f(m)$ be the product of all factors that depend on $m$: \[f(m)=\frac{1}{m} \prod_{\substack{l\mid m\\l \nmid b'\\l\in \mathcal{P}\setminus \mathcal{G}\\2<l \leq x^{1/2}}}\left( 1-\frac{1}{l-1} \right)^{-1}.\]

        Recall from equation \eqref{eq:defM} that every $m\in \mathcal{M} $ is even, has no Golomb prime factors and satisfies $(m,b)=D$. Let 
        \[
        D=2^v D'\qquad \text{with}\qquad (2,D')=1,
        \]
         and let $a=\max\{1,v\}$. Since $f$ is multiplicative,
        \begin{align*}
            \sum_{\substack{m\in \mathcal{M}\\1\leq m \leq x/t}}f(m)&\leq (f(2^{a})+f(2^{a+1})+\ldots) \prod_{\substack{l \nmid b'D\\l\in \mathcal{P} \setminus \mathcal{G}\\2<l \leq x/t} }(1+f(l)+f(l^2)+ \ldots)  \\ 
            &\qquad\prod_{\substack{l\mid D\\l\nmid b',\,l\neq 2\\l \in \mathcal{P} \setminus \mathcal{G}}}(f(l^{\nu_{l}(D)})+f(l^{\nu_{l}(D)+1})+ \ldots)\prod_{\substack{l\mid D\\l\mid b',\,l\neq 2 \\l \in \mathcal{P} \setminus \mathcal{G}}}(f(l^{\nu_{l}(D)}))
        \end{align*}
    For the last product, note that if a prime $l$ divides both $D$ and $b'$, then $\nu_{l}(b)>\nu_{l}(D)$, so the condition $(m,b)=D$ forces $\nu_{l}(m)=\nu_{l}(D)$.
      Evaluating these Euler factors, we find that the right-hand side equals
        \[
           \frac{1}{2^{a-1}} \prod_{\substack{l \nmid b'D\\l\in \mathcal{P} \setminus \mathcal{G}\\2<l \leq x/t} }\left( 1+\frac{1}{l-2} \right)  \prod_{\substack{l\mid D\\l\nmid b',\,l\neq 2 \\l \in \mathcal{P} \setminus \mathcal{G}}} \frac{1}{l^{\nu_{l}(D)-1}(l-2)}\prod_{\substack{l\mid D\\l\mid b',\,l\neq 2 \\l\in \mathcal{P}\setminus \mathcal{G}}}  \frac{1}{l^{\nu_{l}(D)}}.
        \]
        Since $1+\frac{1}{l-2}=\left( 1-\frac{1}{l-1} \right)^{-1}$, we get
        \[
      \sum_{\substack{m\in \mathcal{M}\\1\leq m \leq x/t}}f(m)\ll \frac{1}{D} \prod_{\substack{l\in \mathcal{P}\setminus \mathcal{G}\\2<l\leq x/t}} \left(1-\frac{1}{l-1}\right)^{-1} \prod_{\substack{l\mid D,\,l\nmid b'\\l\in \mathcal{P}\setminus \mathcal{G}}} \left(1-\frac{1}{l}\right)^{-1} \prod_{\substack{l\mid b',\,l\neq 2 \\l\in \mathcal{P}\setminus \mathcal{G}}}\left( 1-\frac{1}{l-1} \right)  
        \]
    Inserting this bound into equation \eqref{eq:errortx}, we obtain
    \begin{align*}
         \sum_{\substack{g\in \mathcal{G}\\t<g\leq x}} B(x;b,h;g)&\ll_A \frac{x(\log_{3}x)^2}{b'D} \prod_{\substack{r\in\mathcal{P}\\ 2<r\le x^{1/2}}}\Big(1-\frac2r\Big)\prod_{\substack{r\mid b'\\r\in\mathcal{P}\\ r>2}}\Big(1-\frac2r\Big)^{-1}\prod_{\substack{q\in\mathcal{G}\\ 3<q\le x^{1/2}}}\Big(1-\frac1{q-2}\Big) \\
         &\quad \prod_{\substack{l\in \mathcal{P}\setminus \mathcal{G}\\2<l\leq x/t}} \left(1-\frac{1}{l-1}\right)^{-1} \prod_{\substack{l\mid D,\,l\nmid b'\\l\in \mathcal{P}\setminus \mathcal{G}}} \left(1-\frac{1}{l}\right)^{-1} \prod_{\substack{l\mid b',\,l\neq 2 \\l\in \mathcal{P}\setminus \mathcal{G}}}\left( 1-\frac{1}{l-1} \right)  .
    \end{align*}
    Since $D$ is Golomb-prime free, every prime factor of $D$ lies in $\mathcal{P}\setminus \mathcal{G}$. Using \[\displaystyle \phi(b)=b'D\prod_{\substack{r|b'\\r\in \mathcal{P}}}\left( 1-\frac{1}{r} \right)\prod_{\substack{r\mid D,\,r\nmid b'\\r\in \mathcal{P}}} \left(1-\frac{1}{r}\right)\] and bounding  $\prod_{r\leq x^{1/2}}\left( 1-\frac{1}{r} \right)$ by Mertens' theorem \eqref{mertens2}, we obtain
    \begin{align*}
        \sum_{\substack{g\in \mathcal{G}\\t<g\leq x}} B(x;b,h;g)&\ll_A \frac{{x(\log_{3}x)^2}}{\phi(b)(\log x)} \prod_{\substack{x/t<l\leq x^{1/2}\\l\in \mathcal{P}\setminus \mathcal{G}}}\left( 1-\frac{1}{l-1} \right)\prod_{\substack{q\in\mathcal{G}\\ 3<q\le x^{1/2}}}\Big(1-\frac2{q-1}\Big) \prod_{\substack{q\mid b'\\q\in \mathcal{G}}}\left( 1-\frac{1}{q-1} \right)^{-1} 
    \end{align*}
    The product over $l$ is at most $1$. Since $1-\frac{2}{q-1}<\left( 1-\frac{1}{q-1} \right)^{2},$ the product over $3<q\leq x^{1/2}$ is $\ll(\log_{2}x)^{-2}$. As shown above, the product over $q\mid b'$ is $\ll \log_{3}x.$ Therefore
    \begin{align*}
        \sum_{\substack{g\in \mathcal{G}\\t<g\leq x}} B(x;b,h;g)&\ll \frac{{x(\log_{3}x)^3}}{\phi(b)(\log x)(\log_{2}x)^2} 
    \end{align*}
        \end{proof}
        \begin{remark}
OpenAI has announced a proof that no Dirichlet $L$-function vanishes in the half-plane $\operatorname{Re}s>7/8$ \cite{OpenAI2026QRH}. Using this, one can show that Theorem~\ref{thm:SW} holds uniformly for $1\leq b\leq x^{1/8-2\varepsilon}$, for any fixed $\varepsilon>0$. The proof follows the same steps, with the following changes.
\begin{enumerate}
    \item[(i)] For the main term, we use the bound $|E(x;q,a)|\ll x^{7/8}\log x$, which follows from this zero-free half-plane (a quasi-Riemann hypothesis), instead of the Bombieri--Vinogradov theorem \eqref{eq:BV}. We then apply Theorem~\ref{thm:FLST} with $z_{0}=\exp\big(\frac{\varepsilon\log x}{\log_{3}x}\big)$.
    \item[(ii)] Lemmas~\ref{lem:errorzt} and~\ref{lem:errortx} hold, with the same bounds, for $t=\frac{x}{b}\exp\big(-\frac{\log x}{(\log_{2}x)^{4/5}}\big)$.
    
\end{enumerate}
\end{remark}

\section{Non-smooth shifted Golomb primes: proof of Theorem~\ref{thm:SWnonsmooth}} \label{sec:SWnonsmooth}
In this section we prove Theorem~\ref{thm:SWnonsmooth}, which we restate here.
\mtheoremthree*
We show that among the Golomb primes $g\leq x$ with
$g\equiv h\pmod b$, at most  $26\%$ have $g-h$ divisible by a prime power larger than $x^{0.99}$.

Recall the definition
\[
{G}_{\rm nsm}(x;b,h;z):=\#\{g\in\mathcal{G}(x;b,h):\ g>h\ \text{and}\ g-h\text{ has a prime power divisor}\ >z\},
\]
where $\mathcal{G}(x;b,h)$ is the set of Golomb primes $g\leq x$ with $g\equiv h\pmod{b}$. We first show that the number of positive integers $\leq x$ that are multiples of some prime power $p^k>x^{0.99}$ with $k\geq 2 $ is $O(x^{0.505})$. For each prime $p$, let $k_{p}$ be the least integer such that $p^{k_{p}}>x^{0.99}$. We begin with the primes $p$ with $k_{p}=2$. These primes are greater than $x^{0.495}$, so
    \[\sum_{x^{0.495}\leq p<x^{0.99}} \frac{x}{p^2}=O\left( \frac{x}{x^{0.495}} \right)=O(x^{0.505}).\]
    For primes $\leq x^{0.495}$, we have
    \[\sum_{p<x^{0.495}} \frac{x}{p^{k_{p}}}\leq \sum_{p<x^{0.495}} \frac{x}{x^{0.99}}=O(x^{0.505}). \]
    So prime powers $p^k$ with $k\geq 2$ contribute $O(x^{0.505})$, which is negligible compared with the bound we are proving. It remains to count the $g$ for which $g-h$ has a prime factor $p>x^{0.99}$. Let
     \[{G}'_{\rm nsm}(x;b,h;x^{0.99}):=\#\{g\in\mathcal{G}(x;b,h):\ g>h\ \text{and}\ g-h\text{ has a prime divisor}\ > x^{0.99}\},\]
     so that ${G}_{\rm nsm}(x;b,h;x^{0.99})\leq{G}'_{\rm nsm}(x;b,h;x^{0.99})+O(x^{0.505})$.
        Let $g$ be counted by ${G}'_{\rm nsm}(x;b,h;x^{0.99})$, and let $p>x^{0.99}$ be a prime dividing $g-h$. Since $b\mid g-h$ and $p>b$, we can write
    \[
    g=pmb+h
    \]
    for a positive integer $m$. We fix $m$ and count the possible $p$. Let
     \[N_{m}=\#\left\{  p:x^{0.99}< p<x/mb, pmb+h \text{ is a Golomb prime} \right\}.\]
      Since $p>x^{0.99}$, we have $mb\leq(x+|h|)/x^{0.99}\leq 2x^{0.01}$. Since $p$ and $g$ are primes, $m$ must satisfy the following conditions.

     \begin{enumerate}
         \item $(mb,h)=1$, since a common prime factor would divide $g$.
         \item $(mb,h-1)$ is Golomb-prime-free, since a Golomb prime $q$ dividing both would give $g\equiv 1\pmod q$.
         \item Since $g$ and $p$ are odd, $2\mid mb$ or $2\mid h$.
        \item Since $g>3$ is a Golomb prime, $g\equiv 2\pmod 3$. If $3\nmid h(h-1)$, then $h\equiv 2\pmod 3$, so $3\mid g-h=pmb$, and hence $3\mid mb$. So $3\mid mbh(h-1)$.
     \end{enumerate}
   We therefore define $\mathcal{M}$ as
    \[
    \mathcal{M}:=\{ m: mb\leq 2 x^{0.01}, 2\mid mbh, 3\mid mbh(h-1), (mb,h)=1,(mb,h-1)\ \text{is Golomb-prime-free}\}.
    \]
   If $N_{m}\neq0$, then $m\in\mathcal{M}$. Hence \begin{equation}\label{eq:G'nsm}
   {G}'_{\rm nsm}(x;b,h;x^{0.99})=\sum_{m\in\mathcal{M}}N_{m}.  \end{equation}

Fix $m\in\mathcal{M}$. Since $g$ is a Golomb prime, $g\not\equiv 0\pmod{r}$ for every prime $r<g$, and $g\not\equiv 1\pmod{q}$ for every Golomb prime $q<g$. We translate these conditions on $g$ into conditions on $p$:
    \begin{itemize}
        \item If $r\nmid mb$ and $r\nmid h$, then $g= pmb+h\not\equiv 0 \pmod{r}$ implies $p\not\equiv-\frac{{h}}{mb}\pmod{r}$.
        \item If $q\nmid mb$ and $q\nmid h-1$, then $g=pmb+h\not\equiv 1 \pmod{q}$ implies $p\not\equiv- \frac{{h-1}}{mb}\pmod{q}$.
    \end{itemize}

      Applying Lemma~\ref{lem:sievebound} with $\displaystyle z_{m}=\frac{(x/mb)^{1/4}}{(\log(x/mb))^{100}}\geq \frac{(x^{0.99}/2)^{1/4}}{(\log x)^{100}}=:z$, we obtain
    \[
        N_{m}\leq (1+o(1)) \frac{{2e^{2\gamma}x}}{mb\log\left(  \frac{x}{mb} \right)} \prod_{\substack{l \leq z\\l \nmid mb\\ l \nmid h\\ l\in \mathcal{P} \setminus \mathcal{G}}}\left( 1-\frac{1}{l-1} \right)   \prod_{\substack{q \leq z\\q \nmid mb\\ q \nmid h(h-1)\\ q\in \mathcal{G}}}\left( 1-\frac{2}{q-1} \right)   \prod_{\substack{q \leq z\\q \nmid mb\\ q \mid h(h-1)\\ q\in \mathcal{G}}}\left( 1-\frac{1}{q-1} \right).
    \]
    Since $mb\leq 2x^{0.01}$, we have $\log(x/mb)\geq(0.99+o(1))\log x$. Using this and rearranging the products, we get
    \begin{align*}
        N_{m}&\leq (1+o(1))\frac{2e^{2\gamma}x}{0.99 mb\log x} \prod_{\substack{2<r\leq z\\r\in\mathcal{P}}}\left(1-\frac {1}{r-1}\right) \prod_{\substack{3<q\leq z\\q\in\mathcal{G}}}\left(1-\frac {1}{q-2}\right) \prod_{\substack{r\mid mb\\2<r\leq z\\r\in\mathcal{P}}}\left(1-\frac{1}{r-1}\right)^{-1} \displaybreak[0]\\
        &\quad \prod_{\substack{l\mid h\\2<l\leq z\\l\in\mathcal{P}\setminus\mathcal{G}}}\left(1-\frac{1}{l-1}\right)^{-1} \prod_{\substack{q\mid h(h-1)\\3<q\le z\\q\in\mathcal{G}}}\left(1-\frac{1}{q-2}\right)^{-1} \prod_{\substack{q\mid mb\\3<q\leq z\\q\in\mathcal{G}}}\left(1-\frac{1}{q-2}\right)^{-1}.
    \end{align*}
    Splitting the products over $r\mid mb$ and $q\mid mb$ into primes dividing $b$ and primes dividing $m$ but not $b$, we obtain
    \begin{align}
        N_{m}&\leq (1+o(1))\frac{2e^{2\gamma}x}{0.99 mb\log x} \prod_{\substack{2<r\leq z\\r\in\mathcal{P}}}\left(1-\frac {1}{r-1}\right) \prod_{\substack{3<q\leq z\\q\in\mathcal{G}}}\left(1-\frac {1}{q-2}\right) \prod_{\substack{r\mid b\\2<r\leq z\\r\in\mathcal{P}}}\left(1-\frac{1}{r-1}\right)^{-1} \notag \displaybreak[0]\\
        &\quad \prod_{\substack{r\mid m\\r \nmid b\\2<r\leq z\\r\in\mathcal{P}}}\left(1-\frac{1}{r-1}\right)^{-1} \prod_{\substack{l\mid h\\2<l\leq z\\l\in\mathcal{P}\setminus\mathcal{G}}}\left(1-\frac{1}{l-1}\right)^{-1} \prod_{\substack{q\mid h(h-1)\\3<q\le z\\q\in\mathcal{G}}}\left(1-\frac{1}{q-2}\right)^{-1} \notag \displaybreak[0]\\
        &\quad \prod_{\substack{q\mid b\\3<q\leq z\\q\in\mathcal{G}}}\left(1-\frac{1}{q-2}\right)^{-1} \prod_{\substack{q\mid m\\q \nmid b\\3<q\leq z\\q\in\mathcal{G}}}\left(1-\frac{1}{q-2}\right)^{-1}. \label{eq:Nmwithoutfm}
    \end{align}
    Let $f(m)$ be the product of the factors on the right-hand side of \eqref{eq:Nmwithoutfm} that depend on $m$:
    \[
    f(m)=\frac{1}{m}\prod_{\substack{r\mid m\\r \nmid b\\2<r\leq z\\r\in\mathcal{P}}}\left(1-\frac{1}{r-1}\right)^{-1} \prod_{\substack{q\mid m\\q \nmid b\\3<q\leq z\\q\in\mathcal{G}}}\left(1-\frac{1}{q-2}\right)^{-1}.
    \]
  We need two lemmas. The first bounds the sum of $f(m)$ over $m\in\mathcal{M}$, and the second shows that a product over the Golomb primes in $(2x^{0.01},z]$ is close to $1$.
 \begin{lemma}\label{lem:fm}
    We have
    \[\sum_{m\in \mathcal{M}}f(m)\leq\prod_{\substack{r\in \mathcal{P}\\r\mid b}}\left( 1-\frac{1}{r} \right)^{-1} \prod_{\substack{l \nmid b\\l \nmid h\\l\in \mathcal{P}\setminus \mathcal{G}\\ 2<l \leq 2 x^{0.01}}} \left( 1- \frac{1}{l-1} \right)^{-1}\prod_{\substack{q \nmid b\\q \nmid h(h-1)\\q\in \mathcal{G}\\3<
q \leq 2 x^{0.01}}} \left( 1- \frac{1}{q-2} \right)^{-1}.\]
\end{lemma}
\begin{lemma}\label{lem:h-1Golomb}
    We have
    \[\prod_{\substack{2x^{0.01}<q\leq z\\q\in \mathcal{G}}}\left( 1-\frac{1}{q-1} \right)^{-1} =1+O\left( \frac{1}{\log_{2}x} \right).\]
\end{lemma}
    We first finish the proof of Theorem~\ref{thm:SWnonsmooth}, assuming these two lemmas.
    Summing \eqref{eq:Nmwithoutfm} over $m\in \mathcal{M}$ and using Lemma~\ref{lem:fm} to bound $\sum_{m\in \mathcal{M}} f(m)$, we obtain
    \begin{align*}
        \sum_{m \in \mathcal{M}} N_{m}&\leq (1+o(1))\frac{2e^{2\gamma}x}{0.99 b\log x} \prod_{\substack{2<r\leq z\\r\in\mathcal{P}}}\left(1-\frac {1}{r-1}\right) \prod_{\substack{3<q\leq z\\q\in\mathcal{G}}}\left(1-\frac {1}{q-2}\right) \prod_{\substack{r\mid b\\2<r\leq z\\r\in\mathcal{P}}}\left(1-\frac{1}{r-1}\right)^{-1} \displaybreak[0]\\
        &\quad \prod_{\substack{l\mid h\\2<l\leq z\\l\in\mathcal{P}\setminus\mathcal{G}}}\left(1-\frac{1}{l-1}\right)^{-1} \prod_{\substack{q\mid h(h-1)\\3<q\le z\\q\in\mathcal{G}}}\left(1-\frac{1}{q-2}\right)^{-1} \prod_{\substack{q\mid b\\3<q\leq z\\q\in\mathcal{G}}}\left(1-\frac{1}{q-2}\right)^{-1} \prod_{\substack{r\in \mathcal{P}\\r\mid b}}\left( 1-\frac{1}{r} \right)^{-1} \displaybreak[0]\\
        &\quad \prod_{\substack{l \nmid b\\l \nmid h\\ l\in \mathcal{P}\setminus \mathcal{G}\\2<l \leq 2 x^{0.01}}} \left( 1- \frac{1}{l-1} \right)^{-1} \prod_{\substack{q \nmid b\\q \nmid h(h-1)\\q\in \mathcal{G}\\3<q \leq 2 x^{0.01}}} \left( 1- \frac{1}{q-2} \right)^{-1}.
    \end{align*}
        Combining terms we get
    \begin{align*}
        \sum_{m \in \mathcal{M}} N_{m}&\leq (1+o(1))\frac{2e^{2\gamma}x}{0.99 \phi(b)\log x} \prod_{\substack{q\leq z\\q\in\mathcal{G}}}\left(1-\frac {1}{q-1}\right) \prod_{\substack{q\mid b\\q\in\mathcal{G}}}\left(1-\frac{1}{q-1}\right)^{-1} \prod_{\substack{2x^{0.01}<r \leq z\\r\in \mathcal{P} }} \left( 1- \frac{1}{r-1} \right) \displaybreak[0]\\
        &\quad \prod_{\substack{2x^{0.01}<q\leq z\\q\in\mathcal{G}}}\left(1-\frac {1}{q-1}\right)^{-1} .
    \end{align*}
    By Theorem~\ref{thm:erdos}, the product over $q\leq z$ is $(1+o(1))/\log_{2}x$. By Mertens' theorem \eqref{mertens2},
    \[
    \prod_{\substack{2x^{0.01}<r\leq z\\r\in\mathcal{P}}}\left(1-\frac{1}{r-1}\right)=(1+o(1))\frac{\log(2x^{0.01})}{\log z}=\frac{0.04}{0.99}+o(1),
    \]
    and by Lemma~\ref{lem:h-1Golomb}
    \[
     \prod_{\substack{2x^{0.01}<q\leq z\\q\in\mathcal{G}}}\left(1-\frac{1}{q-1}\right)^{-1}=1+O\left( \frac{1}{\log_{2}x} \right).
    \]
    We conclude that
    \begin{align*}
        {G}'_{\rm nsm}(x;b,h;x^{0.99})= \sum_{m \in \mathcal{M}} N_{m}&\leq(1+o(1))\frac{0.08e^{2\gamma}x}{0.99^2 \phi(b)\log x\log_{2}x} \prod_{\substack{q\mid b\\q\in\mathcal{G}}}\left(1-\frac{1}{q-1}\right)^{-1}  .
    \end{align*}
  Since $0.08e^{2\gamma}/0.99^{2}<0.259$ and ${G}_{\rm nsm}(x;b,h;x^{0.99})\leq{G}'_{\rm nsm}(x;b,h;x^{0.99})+O(x^{0.505})$, this proves Theorem~\ref{thm:SWnonsmooth}.

\begin{proof}[Proof of Lemma~\ref{lem:fm}]
   Since $f$ is multiplicative, we have \[\sum_{m\in \mathcal{M}}f(m)\leq\prod_{\substack{r\in \mathcal{P}\\r\leq 2 x^{0.01}}}F(r),\]
    where $F(r)$ is defined below for each prime $r$:

    \begin{itemize}
        \item For a prime $l$ that is not a Golomb prime, with $l \not=2$, $l \nmid b$ and $l \nmid h$, \[ F(l)=1+f(l)+f(l^2)+ \ldots=1+ \frac{1}{l-2}=\left( 1- \frac{1}{l-1} \right)^{-1}.\]
        \item For a prime $l$ that is not a Golomb prime, with $l \not=2$, $l\mid b$ and $l \nmid h$, \[F(l)=1+f(l)+f(l^2)+ \ldots=1+\frac{1}{l}+\frac{1}{l^2}+ \ldots=\left( 1-\frac{1}{l} \right)^{-1}.\]
        \item For $l=2$, we have $F(2)=1$  if $2\nmid bh$ (as in this case $2\mid m$ for every $m\in\mathcal{M}$), $F(2)=1$ if $2\mid h$, and $F(2)=2=\left( 1-\frac{1}{2} \right)^{-1}$ if $2\mid b$.

        \item For a Golomb prime $q$ with $q\not=3$, $q \nmid b$ and $q \nmid h(h-1)$,
        \[F(q)=1+f(q)+f(q^2)+ \ldots=1+\frac{1}{q-3}=\left( 1-\frac{1}{q-2} \right)^{-1}.\]
          \item For a Golomb prime $q$ with $q\not=3$, $q\mid b$ and $q \nmid h(h-1)$,
        \[F(q)=1+f(q)+f(q^2)+ \ldots=1+\frac{1}{q}+\frac{1}{q^2}+ \ldots =\left( 1-\frac{1}{q} \right)^{-1}.\]
         \item For $q=3$, we have $F(3)\leq1$ if $3\nmid bh(h-1)$ (in which case $3\mid m$ for every $m\in\mathcal{M}$), $F(3)=1$ if $3\mid h(h-1)$, and $\displaystyle F(3)=\frac{3}{2}=\left( 1-\frac{1}{3} \right)^{-1}$ if $3\mid b$.
          \item If $r\mid h$, or if $r$ is a Golomb prime dividing $h-1$, then $r\nmid m$ for every $m\in\mathcal{M}$, so $F(r)=1$.
    \end{itemize}
    Combining these, we obtain
    \[\sum_{m\in \mathcal{M}}f(m)\leq \prod_{\substack{r\in \mathcal{P}\\r\mid b}}\left( 1-\frac{1}{r} \right)^{-1} \prod_{\substack{l \nmid b\\l \nmid h\\l\in \mathcal{P}\setminus \mathcal{G}\\ 2<l \leq 2 x^{0.01}}} \left( 1- \frac{1}{l-1} \right)^{-1}\prod_{\substack{q \nmid b\\q \nmid h(h-1)\\q\in \mathcal{G}\\3<
    q \leq 2 x^{0.01}}} \left( 1- \frac{1}{q-2} \right)^{-1}.\]
    This finishes the proof of Lemma~\ref{lem:fm}.
    \end{proof}
     \begin{proof}[Proof of Lemma~\ref{lem:h-1Golomb}]
    Taking logarithms, we have
    \[
        \log\Bigg(\prod_{\substack{2x^{0.01}<q\leq z\\q\in \mathcal{G}}}\left( 1-\frac{1}{q-1} \right)^{-1}\Bigg)=\sum_{\substack{2x^{0.01}<q\leq z\\q\in \mathcal{G}}} -\log\left( 1-\frac{1}{q-1} \right).
    \]
    Expanding the logarithm, we see that the right-hand side equals
    \[
        \sum_{\substack{2x^{0.01}<q\leq z\\q\in \mathcal{G}}} \frac{1}{q}+O(x^{-0.01}),
    \]
    which, by Lemma~\ref{lem:partialsum}, is
    \[
        (1+o(1))\big(\log_{3}z-\log_{3}(2x^{0.01})\big)+O\left( \frac{1}{\log x\log_{2}x} \right).
    \]
    Since $\displaystyle z=\frac{x^{0.99/4}}{(\log x)^{100}}$, we have
    \[\log_{3}(z)-\log_{3}(2x^{0.01})=O\left( \frac{1}{\log_{2}x} \right).\]
    
    Substituting these estimates, we obtain
    \[
        \log\Bigg(\prod_{\substack{2x^{0.01}<q\leq z\\q\in \mathcal{G}}}\left( 1-\frac{1}{q-1} \right)^{-1}\Bigg)=O\left( \frac{1}{\log_{2}x} \right),
    \]
    and therefore
    \[
        \prod_{\substack{2x^{0.01}<q\leq z\\q\in \mathcal{G}}}\left( 1-\frac{1}{q-1} \right)^{-1}=1+O\left( \frac{1}{\log_{2}x} \right).
    \]
    This completes the proof of Lemma~\ref{lem:h-1Golomb}.
    \end{proof}

   \section{Proof of Theorem~\ref{thm:main}} \label{mainthm}
   The goal of this section is to prove the main result which we restate below:
    \mtheoremone*
    We first explain why the conditions on $b$ in the above theorem, Theorem~\ref{thm:main}, are necessary. Suppose that $(a,b)=1$ and that \eqref{eq:Golrep} holds for a finite set $B$ of Golomb primes, all greater than $|h|$. Then every prime factor of $b$ divides $g-h$ for some $g\in B$.
\begin{itemize}
    \item If a prime $p$ divides both $b$ and $h$, then $p\mid g-h$ gives $p\mid g$, so $p=g$. This is impossible, since $g>|h|\geq p$. So $(b,h)=1$.
    \item If a Golomb prime $q$ divides both $b$ and $h-1$, then $q\mid g-h$ gives $q\mid g-1$. This is impossible: if $g>q$, it contradicts the definition of the Golomb primes, and if $g\leq q$, then $0<g-1<q$. Hence $(b,h-1)$ is Golomb-prime-free.
\end{itemize}
We also note that it is enough to prove Theorem~\ref{thm:main} when $a/b\leq 1$. Indeed, suppose that $a/b>1$. Subtracting $1$ from $a/b$ repeatedly, we reach a fraction $a'/b$ with $0<a'/b\leq 1$, so that
\[
\frac{a}{b}=1+\cdots+1+\frac{a'}{b}.
\]
Here $a'=a-kb$ for some positive integer $k$, so the denominator is still $b$ and $(a',b)=(a,b)=1$.   Both $a'/b$ and $1$ are at most $1$ and satisfy the hypotheses of Theorem~\ref{thm:main}, so we can apply the case $a/b\leq 1$ to each of them. Each time, we choose $y$ larger than every Golomb prime used so far, so the new Golomb primes are different from the earlier ones. Adding these representations gives $a/b$ as a sum of distinct unit fractions with denominators of the form $g-h$, where every $g$ is greater than the original $y$. So it is enough to prove Theorem~\ref{thm:main} for $a/b\leq 1$, and from now on we assume that $a/b\leq 1$.

Our main tool is the following proposition. It is a variant of \cite[Proposition~5.2]{LS24} and \cite[Proposition~3.1]{Bloom2025density}, both of which go back to \cite[Proposition~3]{Croot2003}. The difference is that we need condition~(iii) with a smaller value of $\eta$ than \cite[Proposition~5.2]{LS24} allows. We obtain this by using the fact that $A$ consists of shifted primes.

In the proposition we use the notation as stated in Section \ref{sec:notation}: for a finite set $A$ of positive integers,
\[
R(A)=\sum_{n\in A} \frac{1}{n}, 
\]
and for every prime power $q$, we write $q\pdiv n$ if $q\mid n$ and $(q,n/q)=1$. Then $A_{q}$ is the set of $n\in A $ with $q\pdiv n$, 
\[
R(A;q)=\sum_{n\in A_{q}} \frac{q}{n}, 
\]
and $\mathcal{Q}_{A}$ is the set of prime powers $q$ with $q\pdiv n$ for some $n\in A$. In this section, $q$ denotes a prime power unless we say that $q\in \mathcal{G}$.
\begin{proposition}\label{prop:main}    Let $h$ be a nonzero integer, $N$ be sufficiently large in terms of $h$ and let $N^{0.997}\leq M\leq N/10^4$. Let $A\subseteq\{p-h:\ p \text{ prime}\}\cap(M,N]$ satisfy the following conditions.
    \begin{enumerate}
        \item[(i)] Every $n\in A$ is $N^{0.99}$-power-smooth, in the sense that every prime power dividing $n$ is at most $N^{0.99}$.
        \item[(ii)] $ \frac {99} {100}\log_2N\leq \omega(n)\leq 2\log_2N$ for all $n\in A$.
        \item[(iii)] $R(A;q)\geq \eta=(\log N)^{-1/200}$ for all $q\in\mathcal{Q}_A$.
        \item[(iv)] $R(A)\geq (\log N)^{-1/101}$.
    \end{enumerate}
    Let $Q=\lcm(A)$. If $x\in\mathbb{N}$ satisfies $\displaystyle \frac{x}{Q}\in\bigg[\frac{R(A)}{\log N},\left( \frac{{\log N}}{\log N+1} \right)R(A)\bigg]$, then there exists $B\subseteq A$ with $R(B)=x/Q$.
\end{proposition}
We prove Proposition~\ref{prop:main} in Appendix~\ref{appendix}. In the rest of this section we deduce Theorem~\ref{thm:main} from it.

Let $N$ be large in terms of $a$, $b$, $h$ and $y$. We use Proposition~\ref{prop:main} to write $a/b$ as a sum of reciprocals of $h$-shifted Golomb primes, that is, positive integers $g-h$, with $g\in \mathcal{G}$ . The rough plan is as follows.
\begin{itemize}
    \item[(i)] We break the $h$-shifted Golomb primes into blocks $(N_{j+1},N_{j}]$ for $1\leq j\leq J$, where $N_{1}=N$, $N_{j+1}=N_{j}^{e^{-0.003}}$ and $J=\lfloor 333\log_{2}N\rfloor$.
    \item[(ii)] The reciprocal sum of the $h$-shifted Golomb primes in each block is $\Theta\left( \frac{1}{\log_{2}N} \right)$.
    \item[(iii)] We prune each block to a set $A_{j}$ satisfying the hypotheses of Proposition~\ref{prop:main}, with $R(A_{j})=\Theta\left( \frac{1}{\log_{2}N} \right)$.
    \item[(iv)] We show that a fixed modulus $P$ divisible by $b$, of size about $(\log_{2}N)^2$, divides $\lcm(A_{j})$ for all $j$.
    \item[(v)] We represent suitable fractions $y_{j}/P$ using the $h$-shifted Golomb primes in $A_{j}$.
    \item[(vi)] Finally, we add the representations of $y_{j}/P $ to obtain a representation of $a/b$.
\end{itemize}
For $1\leq j\leq J$, let $A_{j}'$ be the set of $h$-shifted Golomb primes in the interval $(N_{j+1},N_{j}]$.

   \begin{lemma}\label{lem:RSA'} For fixed $h$, we have
       \[R(A_{j}')=\sum_{\substack{N_{j+1}<g-h\leq N_{j}\\g\in \mathcal{G}}} \frac{1}{g-h}=(1+o(1))\frac{0.003}{\log_{2}N_{j}}, \]
       uniformly for $1\leq j\leq J$.
   \end{lemma}
   \begin{proof}
       We have
       \begin{align*}
           \sum_{\substack{N_{j+1}<g-h\leq N_{j}\\g\in \mathcal{G}}} \frac{1}{g-h}&=  \sum_{\substack{N_{j+1}<g-h\leq N_{j}\\g\in \mathcal{G}}} \frac{1}{g}+O_{h}\bigg ( \sum_{\substack{N_{j+1}<g-h\leq N_{j}\\g\in \mathcal{G}}}\frac{1}{g^2} \bigg) \\
           &= \sum_{\substack{N_{j+1}<g\leq N_{j}\\g\in \mathcal{G}}} \frac{1}{g}+\sum_{\substack{N_{j}<g\leq N_{j}+h\\g\in \mathcal{G}}} \frac{1}{g}-\sum_{\substack{N_{j+1}<g\leq N_{j+1}+h\\g\in \mathcal{G}}} \frac{1}{g}+O_{h}\bigg ( \sum_{\substack{N_{j+1}<g-h\leq N_{j}\\g\in \mathcal{G}}}\frac{1}{g^2} \bigg)  \\
           &=\sum_{\substack{N_{j+1}<g\leq N_{j}\\g\in \mathcal{G}}} \frac{1}{g}+O_{h}(N_{j}^{-0.9}).
       \end{align*}
       Since $N_{j+1}=N_{j}^{e^{-0.003}}$, we have $\log_{2}N_{j+1}=\log_{2}N_{j}-0.003$, and hence
       \[\log_{3}N_{j+1}= \log_{3}N_{j}- \frac{0.003}{\log_{2}N_{j}}+O\left( \frac{1}{(\log_{2}N_{j})^2} \right).\]
       Therefore, by Lemma~\ref{lem:partialsum},
       \begin{align*}
              \sum_{\substack{N_{j+1}<g-h\leq N_{j}\\g\in \mathcal{G}}} \frac{1}{g-h}&=(1+o(1))(\log_{3}N_{j}-\log_{3}N_{j+1})+O\left( \frac{1}{\log N_{j+1}\log_{2}N_{j+1}} \right)+O_{h}(N_{j}^{-0.9}) \\
              &= (1+o(1))\frac{0.003}{\log_{2}N_{j}}.\qedhere
       \end{align*}
   \end{proof}
 We now prove three estimates, Lemmas~\ref{lem:w0.99}, \ref{lem:w2} and~\ref{lem:qeta}, which we use to show that the pruning steps below remove only a small part of the reciprocal sum.
   \begin{lemma}\label{lem:w0.99} There is a constant  $\delta>0.00005$ such that
       \[\#\{p\leq N:\ \omega(p-h)<0.99\log_{2}N\}\ll_{h} \frac{{N(\log_{2}N)^{O(1)}}}{(\log N)^{1+\delta}}.\]
   \end{lemma}
   \begin{proof}
       By \cite[Corollary~1]{Ford2022}, for every nonzero integer $a$ and every $k\geq 1 $,
       \[\#\left\{  -a<p\leq x: \omega\left( \frac{{p+a}}{s}\right)=k \right\}\ll \frac{|a|}{\phi(|a|)}\frac{x(\log \log x+O(1))^{k-1}}{(k-1)!(\log x)^2},\]
       where $s=(2,a+1)$. We take $a=-h$ so that $s=(2,h+1)$. If $h$ is even then $s=1$. If $h$ is odd, then $s=2$ and $\omega\left( \frac{{p-h}}{2} \right)\leq\omega(p-h)$. In both cases  $\omega(p-h)<0.99\log_{2}N$ implies  $\omega\left( \frac{{p-h}}{s} \right)<0.99\log_{2}N$.
      Summing over $k<0.99\log_{2}N$, we obtain
       \[
           \#\{ h<p\le N:\ \omega(p-h)<0.99\log_{2}N\}\ll_{h} \sum_{k<0.99\log_{2}N} \frac{N(\log \log N+O(1))^{k-1}}{(k-1)!(\log N)^2}.
       \]
       The term with $k=1$ is $N/(\log N)^2,$ which is smaller than the bound we want. For the remaining terms, we use the approximation $\log(t!)=t\log t-t+O(\log t)$ and see that their sum is
       \begin{align*}
         &\ll \frac{N}{(\log N)^2} \sum_{2\leq k<0.99\log_{2}N}\exp\left( (k-1)\left( \log_{3}N-\log{(k-1)}+1+O\left( \frac{1}{\log_{2}N}+ \frac{\log(k-1)}{k-1} \right) \right) \right) \\
         &\ll \frac{N(\log_{2}N)^{O(1)}}{(\log N)^2} \sum_{2\leq k<0.99\log_{2}N}\exp\left( (k-1)\left( \log_{3}N-\log{(k-1)}+1 \right) \right)  
         \end{align*}
         The function $f(t)=t(\log_{3}N-\log t+1)$ is increasing for $0<t<\log_{2}N$, so each term in the sum is at most $\exp(0.99(1-\log{0.99})\log_{2}N)$. Since there are at most $\log_{2}N$ terms, the sum is 
         \begin{align*}
         &\ll \frac{N(\log_{2}N)^{O(1)}}{(\log N)^2} \exp(0.99(1-\log{0.99})\log_{2}N)\\
         &=\frac{{N(\log_{2}N)^{O(1)}}}{(\log N)^{1+\delta}},
       \end{align*}
       where $\delta>0.00005$.
   \end{proof}
   \begin{lemma}\label{lem:w2} We have
        \[\#\{p\leq N:\ \omega(p-h)>2\log_{2}N\}\ll_{h} \frac{N(\log_{2}N)^{O(1)}}{(\log N)^{2\log2}}.\]
   \end{lemma}
   \begin{proof}
       As in the proof of Lemma~\ref{lem:w0.99}, we apply \cite[Corollary 1]{Ford2022} with $a=-h$, so that $s=(2,h+1)$. Since $s\leq 2 $, we have $\omega\left( \frac{{p-h}}{s} \right)\geq\omega(p-h)-1$, so every prime $p>h$ with $\omega(p-h)>2\log_{2}N$ satisfies $\omega\left( \frac{{p-h}}{s} \right)>2\log_{2}N-1$, hence
       \[
           \#\{ h<p\le N:\ \omega(p-h)>2\log_{2}N\}\ll_{h} \sum_{k>2\log_{2}N-1} \frac{N(\log \log N+O(1))^{k-1}}{(k-1)!(\log N)^2}.
       \]
       
       Let $k_{0}$ be the least integer such that $k_{0}>2\log_{2}N-1$. For $k\geq k_{0},$ the ratio of consecutive terms in this sum is at most
       \[
        \frac{{\log_{2}N+O(1)}}{k-1}\leq \frac{{\log_{2}N+O(1)}}{2\log_{2}N}\leq 0.6,
       \]
       for $N$ large. Hence
          \begin{align*}
           \#\{ h<p\le N:\ \omega(p-h)>2\log_{2}N\}&\ll_{h} \frac{N}{(\log N)^2} \frac{(\log_{2} N+O(1))^{k_{0}-1}}{(k_{0}-1)!} (1+0.6+0.6^2+\cdots) \\
           &\ll \frac{N}{(\log N)^2}\frac{(\log_{2}N+O(1))^{k_{0}-1}}{(k_{0}-1)!}.
       \end{align*}
      Note that $k_{0}=2\log_{2}N+O(1)$. Using the approximation $\log(t!)=t\log t-t+O(\log t)$, we see that the right-hand side is
       \begin{align*}
       &\ll\frac{N}{(\log N)^2} \exp\big(2\log_{2}N(\log_{3}N-\log(2\log_{2}N)+1)+O(\log_{3}N)\big) \\
       &\ll\frac{N(\log_{2}N)^{O(1)}}{(\log N)^2} \exp\big(2\log_{2}N(1-\log2)\big) \\
       &=\frac{N(\log_{2}N)^{O(1)}}{(\log N)^2}(\log N)^{2-2\log2} \\
       &=\frac{N(\log_{2}N)^{O(1)}}{(\log N)^{2\log2}}.\qedhere
       \end{align*}
   \end{proof}
   \begin{lemma} \label{lem:qeta}
       Let $N$ be sufficiently large, let $A\subseteq[1,N]$ be a set of positive integers, and let $\eta>0$. There is a subset $A'\subseteq A$ such that\[R(A')\geq R(A)-1.5\eta \log_{2}N\qquad\text{and}\qquad R(A';q)\geq \eta\quad\text{for all } q\in \mathcal{Q}_{A'}.\]
   \end{lemma}

   \begin{proof}
    We construct $A'$ iteratively. Let $A_{0}=A$. If $A_{i}$ has been defined and there is a prime power $q_{i}\in\mathcal{Q}_{A_{i}}$ with $R(A_{i};q_{i})<\eta$, we remove the elements of $(A_{i})_{q_{i}}$, that is, we set $A_{i+1}=A_{i}\setminus (A_{i})_{q_{i}}$. Then
    \[
    R(A_{i+1})=R(A_{i})-\frac{R(A_{i};q_{i})}{q_{i}}\geq R(A_{i})-\frac{\eta}{q_{i}}.
    \]
    Since $A$ is finite and each step removes at least one element, the process stops after finitely many steps, say at step $k$, and we set $A':=A_{k}$. By construction, $R(A';q)\geq\eta$ for all $q\in\mathcal{Q}_{A'}$. Moreover, $q_{i}\notin\mathcal{Q}_{A_{i'}}$ for $i'>i$, so the prime powers $q_{0},\dots,q_{k-1}$ are distinct, and they are at most $N$. By Mertens' theorem, equation~\eqref{mertens1},
    \[
    R(A')\geq R(A)-\sum_{i=0}^{k-1} \frac{\eta}{q_{i}}\geq R(A)-\eta\sum_{\substack{q\leq N\\ q\text{ prime power}}}\frac{1}{q}\geq R(A)-1.5\eta \log_{2}N.\qedhere
    \]
   \end{proof}
   We now prune $A_{j}'$ to obtain a set $A_{j}$ satisfying the hypotheses of Proposition~\ref{prop:main} with $N=N_{j}$ and with $M=N_{j+1}$.
   \begin{enumerate}
       \item We first remove the elements of $A_{j}'$ that are not $N_{j}^{0.99}$-power-smooth, and we denote the set of removed elements by $(A_{j}')_{\rm nsm}$. If $n=g-h\in(A_{j}')_{\rm nsm}$ and $n\leq t\leq N_{j}$, then $g-h$ has a prime power divisor greater than $N_{j}^{0.99}\geq t^{0.99}$. Hence, by Theorem~\ref{thm:SWnonsmooth} with $b=1$ and partial summation,
       \[R((A_{j}')_{\rm nsm})\leq \frac{0.259}{\log N_{j}\log_{2}N_{j}}+\int_{N_{j+1}}^{N_{j}} \frac{0.259}{t\log t\log_{2}t}\,dt\leq(1+o(1)) \frac{0.000777}{\log_{2}N_{j}}.\]
       \item We then remove the elements $n$ with $\omega(n)<0.99\log_{2}N_{j}$ or $\omega(n)>2\log_{2}N_{j}$, and we denote the set of removed elements by $(A_{j}')_{\omega}$. By Lemmas~\ref{lem:w0.99} and~\ref{lem:w2} and partial summation,
       \[
       R((A_{j}')_{\omega})\ll \frac{(\log_{2}N_{j})^{O(1)}}{(\log N_{j})^{\delta}}+\frac{(\log_{2}N_{j})^{O(1)}}{(\log N_{j})^{2\log2-1}}=o\left( \frac{1}{\log_{2}N_{j}} \right).
       \]
       \item Finally, we apply Lemma~\ref{lem:qeta} to the remaining set, with $N_{j}$ in place of $N$ and with $\eta=(\log N_{j})^{-1/200}$, and we let $A_{j}$ be the resulting set. This step decreases the reciprocal sum by at most \[1.5(\log N_{j})^{-1/200}\log_{2}N_{j}=o(1/\log_{2}N_{j}).\]
   \end{enumerate}
   Collecting these bounds and using Lemma~\ref{lem:RSA'}, we get
   \[R(A_{j})\geq(1+o(1)) \frac{0.003}{\log_{2}N_{j}}-(1+o(1)) \frac{0.000777}{\log_{2}N_{j}}-o\left( \frac{1}{\log_{2}N_{j}}\right)=(1+o(1)) \frac{0.002223}{\log_{2}N_{j}}, \]
   uniformly for $1\leq j\leq J$. By construction, $A_{j}$ satisfies conditions (i)--(iii) of Proposition~\ref{prop:main} with $N_{j}$ in place of $N$, and condition (iv) holds since $R(A_{j})\gg 1/\log_{2}N_{j}$.

   Let $x_0=2 \log_{3}N$ and define
\begin{equation*}
P:=b\prod_{\substack{l\leq x_0\\ l\notin\mathcal{G}\\ l\nmid bh}}l,
\end{equation*}
where $l$ runs over the primes $l\leq x_{0}$ that are not Golomb primes and do not divide $bh$. By the prime number theorem and Theorem~\ref{thm:erdos}, the product of the primes $l\leq x_{0}$ that are not Golomb primes is $\exp\big(2\log_{3}N-(2+o(1))\log_{3}N/\log_{5}N\big)$. Leaving out the primes that divide $bh$ changes this by a factor of at most $b|h|$, so for $N$ large,
\begin{equation}\label{eq:P}
(\log_{2}N)^{1.9}\leq P\leq b(\log_{2}N)^{2.1}.
\end{equation}
We now show that $P$ divides $\lcm(A_{j})$ for all $1\leq j\leq J$.
\begin{lemma}\label{lem:P}
    For $N$ sufficiently large in terms of $a$, $b$ and $h$, we have $P\mid Q_{j}$ for every $1\leq j\leq J$, where $Q_{j}=\lcm(A_{j})$.
\end{lemma}
\begin{proof}
    Fix $j$, and let $m$ be either $b$ or one of the primes $l$ in the definition of $P$. Then $(h,m)=1$ and $(h-1,m)$ is Golomb-prime-free. The elements of $A_{j}'$ divisible by $m$ are the numbers $g-h$ with $g\equiv h\pmod{m}$. By Theorem~\ref{thm:SW} and partial summation, their reciprocal sum is
    \begin{equation} \label{recsum}
         (1+o(1))\frac{0.003}{\phi(m)\log_{2}N_{j}}\prod_{\substack{q\mid m\\q\in \mathcal{G}}}\left( 1- \frac{1}{q-1} \right)^{-1},
    \end{equation}
    and by Theorem~\ref{thm:SWnonsmooth} (with $m$ in place of $b$) and partial summation, the non-smooth elements $(A_{j}')_{\rm nsm}$  contribute at most $0.259+o(1)$ times the expression in equation~\eqref{recsum}. On the other hand, the second and third pruning steps remove elements with total reciprocal sum
    \[
    \ll \frac{(\log_{2}N_{j})^{O(1)}}{(\log N_{j})^{\delta}}+\frac{(\log_{2}N_{j})^{O(1)}}{(\log N_{j})^{2\log2-1}}+\frac{\log_{2}N_{j}}{(\log N_{j})^{1/200}}\ll(\log N_{j})^{-\delta/2},
    \]
   which is much smaller than the reciprocal sum of the elements of $A_{j}'$ divisible by $m$. So the elements of $A_{j}$ that are divisible by $m$ have reciprocal sum at least
    \[
    (0.74+o(1))\frac{0.003}{\phi(m)\log_{2}N_{j}}>0.
    \]
    In particular, some element of $A_{j}$ is divisible by $m$, and so $m\mid Q_{j}$. Since $b$ and the primes $l$ in the definition of $P$ are pairwise coprime, it follows that $P\mid Q_{j}$.
\end{proof}

Define $y_{j}:=\lfloor PR(A_{j}) \rfloor -1$. Then
\[
\frac{y_{j}}{P}=R(A_{j})-\varepsilon_{j}\qquad\text{for some } \varepsilon_{j}\in\left[\frac{1}{P},\frac{2}{P}\right).
\]
We now use Proposition~\ref{prop:main}, with $N_{j}$ in place of $N$, to represent $y_{j}/P$ by elements of $A_{j}$. Since $P\mid Q_{j}$ by Lemma~\ref{lem:P}, we have $y_{j}/P=x_{j}/Q_{j}$ with $x_{j}:=y_{j}Q_{j}/P\in\mathbb{N}$, so it remains to check that
 $y_{j}/P$ lies in $\displaystyle\bigg[\frac{R(A_{j})}{\log N_{j}},\left( \frac{{\log N_{j}}}{\log N_{j}+1} \right)R(A_{j})\bigg]$. The lower bound holds since
\[\frac{y_{j}}{P}=R(A_{j})-\varepsilon_{j}\geq R(A_{j})- \frac{2}{P}>\frac{R(A_{j})}{2},\]
as $1/P\leq (\log_{2}N)^{-1.9}$ by equation~\eqref{eq:P} and $R(A_{j})\gg 1/\log_{2}N$. For the upper bound, it is enough to show that $\displaystyle \varepsilon_{j}\geq \frac{R(A_{j})}{\log N_{j}+1}$, and indeed
\[
\varepsilon_{j}\geq\frac{1}{P}\geq \frac{1}{b(\log_{2}N)^{2.1}}\geq \frac{R(A_{j})}{\log N_{j}+1}.
\]

Since $\log_{2}N_{j}=\log_{2}N-0.003(j-1)$ and $J=\lfloor 333\log_{2}N\rfloor$, the lower bound for $R(A_{j})$ gives
\[
\sum_{1\leq j\leq J}R(A_{j})\geq(1+o(1))\sum_{1\leq j\leq J}\frac{0.002223}{\log_{2}N-0.003(j-1)}
=(0.741+o(1))\log 1000>5.
\]
Since $\sum_{j\leq J}\varepsilon_{j}\leq 2J/P=o(1)$ by \eqref{eq:P}, it follows that $\sum_{j\leq J}y_{j}/P> 1 \geq a/b$.
If \[\frac{a}{b}=\sum_{j\leq k} \frac{y_{j}}{P}\] for some $k$, we are done. Otherwise, let $0\leq k<J$ be such that
\[
\sum_{1\leq j\leq k} \frac{y_{j}}{P}< \frac{a}{b} <\sum_{1\leq j\leq k+1} \frac{y_{j}}{P}.
\]
Since $b\mid P$, the number $\displaystyle\frac{a}{b}-\sum_{1\leq j\leq k}\frac{y_{j}}{P}$ is a positive multiple of $1/P$, so
\[\frac{1}{P}\leq\frac{a}{b}-\sum_{1\leq j\leq k} \frac{y_{j}}{P}< \frac{y_{k+1}}{P}\leq \left( \frac{{\log N_{k+1}}}{\log N_{k+1}+1} \right)R(A_{k+1}).\]
Moreover, since $b\leq \log_{2}N$,
\[
\frac{1}{P}\geq \frac{1}{b(\log_{2}N)^{2.1}}\geq \frac{1}{(\log_{2}N)^{3.1}}\geq \frac{R(A_{k+1})}{\log N_{k+1}}.
\]
Hence, by Proposition~\ref{prop:main} and Lemma~\ref{lem:P}, we can use $A_{k+1}$ to represent $\displaystyle \frac{a}{b}-\sum_{1\leq j\leq k} \frac{y_{j}}{P}$. Since the sets $A_{j}$ are pairwise disjoint, this gives a representation of $a/b$ as a sum of $1/(g-h)$ over distinct Golomb primes $g$. Every element of every $A_{j}$ is greater than $N_{J+1}$, and $N_{J+1}>y+|h|$ once $N$ is large, so all of these Golomb primes are greater than $y$. This completes the proof of Theorem~\ref{thm:main}.
\appendix
\section{Proof of Proposition~\ref{prop:main}}\label{appendix}
The setup follows the proof of \cite[Proposition~5.2]{LS24}, and our treatment of the minor arcs follows \cite[Proposition~3.1]{Bloom2025density}. Throughout this appendix, $q$ denotes a prime power (rather than a Golomb prime, as in other parts of the paper).
Let
\[
\tau=\frac{x/Q}{R(A)}.
\]
Then $\tau\in \displaystyle \left[ \frac{1}{\log N}, \frac{{\log N}}{1+\log N} \right]$ by the assumption on $x/Q$.
Let $B\subseteq A$ be a random subset, where each $n\in A$ is included in $B$ independently with probability $\tau$. Then $\mathbb{E}(R(B))=\tau R(A)=x/Q$. By orthogonality,
\[
 \mathds{1}_{R(B)-x/Q\in\mathbb{Z}}=\frac 1Q\sum_{-Q/2<t\leq Q/2}e\bigg(t\left(R(B)-\frac xQ\right)\bigg).
\]
Taking expectations and using independence, we get
\begin{align*}
    \mathbb{P}\left(R(B)-\frac xQ\in\mathbb{Z}\right)&=\frac 1Q\sum_{-Q/2<t\leq Q/2} e\left(-\frac{tx}{Q}\right) \prod_{n\in A} \mathbb{E}\,e\left(\frac{t\mathds{1}_{n\in B}}{n}\right)\\
    &= \frac 1Q\sum_{-Q/2<t\leq Q/2} e\left(-\frac{tx}{Q}\right) \prod_{n\in A} \bigg(1-\tau+\tau e\left(\frac tn\right)\bigg).
\end{align*}
Since the left-hand side is real, we may take real parts on the right-hand side, and therefore
\begin{equation}\label{eq:Pz}
    \mathbb{P}\left(R(B)-\frac xQ\in\mathbb{Z}\right)= \frac 1Q\sum_{-Q/2<t\leq Q/2} \operatorname{Re}\left( e\left(-\frac{tx}{Q}\right) \prod_{n\in A} \bigg(1-\tau+\tau e\left(\frac tn\right)\bigg)\right).
\end{equation}
We will prove later, using a major-arc and minor-arc analysis, that the right-hand side of \eqref{eq:Pz} is at least $1/(2Q)$. For now, we show that $\mathbb{P}(|R(B)-x/Q|\geq 1)<\frac{1}{4Q}$, as together with the lower bound it implies
\[
\mathbb{P}(R(B)=x/Q)>\frac{1}{2Q}-\frac{1}{4Q}=\frac{1}{4Q}>0.
\]
By condition (i) of Proposition~\ref{prop:main}, every $q\in\mathcal{Q}_{A}$ is at most $N^{0.99}$, and hence
\begin{equation}\label{eq:boundQ}
    Q\leq \prod_{\substack{q\leq N^{0.99}\\q\text{ prime power}}} q\leq e^{2N^{0.99}}.
\end{equation}
Using Hoeffding's inequality \cite[Theorem~2]{Hoeffding1963} for the independent random variables $\mathds{1}_{n\in B}/n\in[0,1/n]$, together with $\mathbb{E}(R(B))=x/Q$ and \eqref{eq:boundQ}, we get
\[
\mathbb{P}(|R(B)-x/Q|\geq 1)\leq 2 \exp\left( -\frac{2}{\sum_{n\in A}1/n^2 } \right)\leq 2 \exp\left( -\frac{2M^2}{N } \right) \leq 2 \exp(-2N^{0.994})<\frac{1}{4Q}.
\]
It remains to show that
\begin{equation}\label{eq:goal}
    \sum_{-Q/2<t\leq Q/2}\operatorname{Re}\left( e\left(-\frac{tx}{Q}\right) \prod_{n\in A} \bigg(1-\tau+\tau e\left(\frac tn\right)\bigg)\right)\geq \frac{1}{2}.
\end{equation}
We do so by splitting the range of $t$ into the major arc and the minor arc,
\[\mathfrak{M}=\big\{t\in(-\tfrac Q2,\tfrac Q2]:|t|\leq \tfrac M2\big\}, \qquad \mathfrak{m}=\big\{t\in(-\tfrac Q2,\tfrac Q2]:|t|> \tfrac M2\big\}.\]
We will use the following standard estimates.
\begin{fact}[{\cite[Fact~2.5]{LS24}}]\label{lem:expestimates}
    Let $\theta\in[-1/2, 1/2]$ and $\lambda\in[0,1]$. Then
    \begin{enumerate}
        \item $\big|(1-\lambda)+\lambda e(\theta)-e(\lambda\theta)\big(1-2\pi^{2} \lambda(1-\lambda)\theta^2\big)\big|\ll |\theta|^3$;
        \item $|1-\lambda+\lambda e(\theta)|\leq 1-8\lambda(1-\lambda)\theta^2\leq \exp(-8\lambda(1-\lambda)\theta^2)$.
    \end{enumerate}
\end{fact}
Since $A\subseteq(M,N]$, condition (iv) gives $|A|\geq MR(A)\geq N^{0.997}(\log N)^{-1/101}\geq N^{0.95}$. Moreover, $\tau\in[(\log N)^{-2},1-(\log N)^{-2}]$. By~\cite[Lemma~3.1]{LS24}, whose proof uses Fact~\ref{lem:expestimates}, we have
\begin{equation}\label{eq:majorgoal}
    \sum_{t\in\mathfrak{M}}\operatorname{Re}\left( e\left(-\frac{tx}{Q}\right) \prod_{n\in A} \bigg(1-\tau+\tau e\left(\frac tn\right)\bigg)\right)\geq \frac{3}{4}.
\end{equation}

We now turn to the minor arc, where we closely follow Bloom~\cite{Bloom2025density}. For each $t\in\mathfrak{m}$, let \[I_{t}=\left\{  m\in\mathbb{Z}:|m-t|< \frac{K}{2}  \right\}, \quad \text{where } K=MN^{-2/\log_{2}N}.\]
Define \[A_{I_{t}}=\left\{  n\in A:\ \|t/n\|<\frac{K}{2n}  \right\}=\{ n\in A:\ n \mid m  \text{ for some } m\in I_{t}\},\]
and
\[\mathfrak{m}_{1}=\left\{  t\in \mathfrak{m}:\ |A\setminus A_{I_{t}}|\geq \frac{M}{\log N}  \right\},\qquad \mathfrak{m}_{2}=\mathfrak{m}\setminus \mathfrak{m}_{1}.\]
For ease of notation, we define
\[F(t)=e\left(-\frac{tx}{Q}\right) \prod_{n\in A} \bigg(1-\tau+\tau e\left(\frac tn\right)\bigg).\]
Since $\operatorname{Re}F(t)\geq -|F(t)|$, by \eqref{eq:goal} and \eqref{eq:majorgoal} it suffices to show that $\sum_{t\in\mathfrak{m}}|F(t)|\leq 1/4$. By Fact~\ref{lem:expestimates}(2), we have $|1-\tau+\tau e(t/n)|\leq \exp(-8\tau(1-\tau)\|t/n\|^{2})$ for all $n\in A$. If $t\in\mathfrak{m}_{1}$, then $\|t/n\|\geq K/(2n)\geq K/(2N)$ for at least $M/\log N$ elements $n\in A$. Hence, using $Q\leq e ^{2N^{0.99}}$ in the last step, we have
\begin{align*}
    \sum_{t\in\mathfrak{m}_{1}}|F(t)|&\leq \sum_{t\in\mathfrak{m}_{1}} \exp\left( -8\tau(1-\tau)\sum_{n\in A}\|t/n\|^2  \right) \leq  \sum_{t\in\mathfrak{m}_{1}}\exp\left( -8\tau(1-\tau) \frac{K^2M}{4N^2\log N} \right)  \\
    &\leq \sum_{-Q/2<t\leq Q/2}e^{-2N^{0.9905}}\leq \frac{1}{8}.
\end{align*}
It remains to show that $\sum_{t\in\mathfrak{m}_{2}} |F(t)|\leq \frac18$.

We say that $q\in \mathcal{Q}_{A}$ is $t$-poor if
\[
\#\left\{  n\in A_{q}:\ \|t/n\|\geq \frac{K}{2n}  \right\}< \frac{{25N^2\log N \log_{2}N}}{\tau(1-\tau)K^2}=:L,
\]
and we let $\mathcal{D}_{t}\subseteq \mathcal{Q}_{A}$ be the set of $t$-poor prime powers.
  If $q\in\mathcal{Q}_A\setminus\mathcal{D}_t$, then at least $L$ elements $n\in A_q$ satisfy $\|t/n\|\geq K/(2n)\geq K/(2N)$, so by Fact~\ref{lem:expestimates}(2),
  \[
  \prod_{n\in A_q}\big|1-\tau+\tau e(t/n)\big|\leq\exp\Big(-\frac{2\tau(1-\tau)LK^{2}}{N^{2}}\Big)=N^{-50\log_2N}.
  \]

Since each factor $|1-\tau+\tau e(t/n)|$ is at most $1$, and every $n\in A$ lies in exactly $\omega(n)\leq 2\log_{2}N$ of the sets $A_{q}$ (recall that $n\in A_{q}$ exactly when $q\pdiv n$), we have
\begin{align*}
    |F(t)|^{2\log_2{N}}&\leq\prod_{n\in A}\big|1-\tau+\tau e(t/n)\big|^{\omega(n)}= \prod_{q\in\mathcal{Q}_A}\ \prod_{n\in A_{q}}\big|1-\tau+\tau e(t/n)\big|\\
    &\leq \prod_{q\in\mathcal{Q}_A\setminus\mathcal{D}_t}N^{-50\log_2{N}}=N^{-50\log_2{N}\,|\mathcal{Q}_A\setminus \mathcal{D}_t|}.
\end{align*}
This gives
\begin{equation}\label{eq:minorbound}
    |F(t)|\leq N^{-25|\mathcal{Q}_A\setminus \mathcal{D}_t|}\qquad \text{for }t\in\mathfrak{m}_{2}.
\end{equation}
Recall that $[\mathcal{D}]$ denotes the lcm of a finite set $\mathcal{D}$ of integers.
\begin{lemma}\label{lem:divisibility}
    For every $t\in\mathfrak{m}_2$, there exists $x(t)\in I_t$ such that $q\mid x(t)$ for all $q\in\mathcal{D}_t$; equivalently, $[\mathcal{D}_t]\mid x(t)$.
\end{lemma}
We first complete the proof of Proposition~\ref{prop:main} assuming Lemma~\ref{lem:divisibility}, and then prove the lemma.
Let $t\in\mathfrak{m}_{2}$. We first note that $\mathcal{D}_{t}\neq\mathcal{Q}_{A}$. Indeed, otherwise $Q=[\mathcal{Q}_{A}]$ divides $x(t)$.
However, $x(t)\neq0$, since $|t|>M/2>K/2$, and $|x(t)|<|t|+K/2\leq Q/2+K/2<Q$. This contradicts $Q\mid x(t)$.

Now let $1\leq s\leq |\mathcal{Q}_{A}|$, and let us count the $t\in\mathfrak{m}_{2}$ with $|\mathcal{Q}_{A} \setminus \mathcal{D}_{t}|=s$. Since $|\mathcal{Q}_{A}|<N$, there are fewer than $N^{s}$ possibilities for the set $\mathcal{Q}_{A}\setminus\mathcal{D}_{t}$, and hence for $\mathcal{D}_{t}$. Fix such a set $\mathcal{D}$. By Lemma~\ref{lem:divisibility}, every $t\in\mathfrak{m}_{2}$ with $\mathcal{D}_{t}=\mathcal{D}$ lies within distance $K/2$ of a multiple of $[\mathcal{D}]$ in $[-Q,Q]$. There are at most $2Q/[\mathcal{D}]+1\leq 3Q/[\mathcal{D}]$ such multiples, and each of them accounts for at most $K+1$ values of $t$. Since $Q$ divides $[\mathcal{D}]\prod_{q\in\mathcal{Q}_{A}\setminus\mathcal{D}}q$ and every $q\in\mathcal{Q}_{A}$ is at most $N^{0.99}$, the number of such $t$ is at most
\[
\frac{3Q}{[\mathcal{D}]}(K+1) \leq 3N^{0.99s}(K+1)\leq N^{s+1}.
\]
Using \eqref{eq:minorbound}, we get
\[
 \sum_{t\in\mathfrak{m}_{2}} |F(t)|\leq  \sum_{t\in\mathfrak{m}_{2}}N^{-25|\mathcal{Q}_A\setminus \mathcal{D}_t|} \leq \sum_{s=1}^{N} N^{-25s} N^{s+1} N^s\leq \sum_{s=1}^{N} N^{-22s}\leq N^{-21}\leq \frac18.
\]
Together with the bound for $\mathfrak{m}_{1}$, this gives $\sum_{t\in\mathfrak{m}}|F(t)|\leq 1 /4$, which proves \eqref{eq:goal}. This completes the proof of Proposition~\ref{prop:main}, assuming Lemma~\ref{lem:divisibility}.

\subsection{Proof of Lemma~\ref{lem:divisibility}}
We will use the following technical lemmas. The first three are Lemmas~4.1, 4.2 and 4.3 of Bloom~\cite{Bloom2025density}, and the fourth follows from a sieve estimate, Lemma~\ref{lem:sievebound}.

\begin{lemma}\label{lem-basic}
If $0<\lvert n_1-n_2\rvert\leq N$, then
\[
  \sum_{q\mid (n_1,n_2)}\frac{1}{q}\ll \log_{3}N,
\]
where the summation is restricted to prime powers.
\end{lemma}

\begin{lemma}\label{lem-rtop}
Let $\frac{1}{2}>\varepsilon>0$ and let $N$ be sufficiently large, depending
on $\varepsilon$. If $A$ is a finite set of integers such that
$R(A)\geq (\log N)^{-\varepsilon/2}$ and
$(1-\varepsilon)\log_{2}N\leq \omega(n)\leq 2\log_{2}N$ for all $n\in A$, then
\[
  \sum_{q\in\mathcal{Q}_A}\frac{1}{q}\geq (1-2\varepsilon)e^{-1}\log_{2}N.
\]
\end{lemma}

\begin{lemma}\label{lem-usingq}
There is a constant $c>0$ such that the following holds. Let
$N\geq M\geq N^{1/2}$ be sufficiently large, and suppose that
$1\leq k\leq c\log_{2}N$. Suppose that $A\subset[M,N]$ is a set of integers
such that $\omega(n)\leq (\log N)^{1/k}$ for all $n\in A$.

For all $q$ such that $R(A;q)\geq (\log N)^{-1/2}$, there exists $d$ such that
\begin{enumerate}
  \item[(i)] $qd>M\exp(-(\log N)^{1-1/k})$,
  \item[(ii)] $\omega(d)\leq \frac{5}{\log k}\log_{2}N$, and
  \item[(iii)] $\sum\limits_{\substack{n\in A_q,\ qd\mid n\\ (qd,n/qd)=1}}\frac{qd}{n}
         \gg \frac{R(A;q)}{(\log N)^{2/k}}$.
\end{enumerate}
\end{lemma}

\begin{lemma}\label{lem:shiftedsieve}
There is an absolute constant $C_{2}>0$ such that the following holds. Let $h\in\mathbb{Z}\setminus\{0\}$. For all $X$ sufficiently large in terms of $h$ and every set $\mathcal{P}'$ of primes $p$ with $3\leq p\leq X^{1/10}$ and $p\nmid h$, we have
\[
\#\big\{n\in[X,2X):\ n+h \text{ is prime and } p\nmid n \text{ for all } p\in\mathcal{P}'\big\}\leq C_2\,\frac{X}{\log X}\prod_{p\in\mathcal{P}'}\Big(1-\frac{1}{p-1}\Big).
\]
\end{lemma}
  This follows from Lemma~\ref{lem:sievebound}, applied to the primes up to $2X+h$, with the single residue class $h\pmod p$ removed for each $p\in\mathcal{P}'$; this class is coprime to $p$ because $p\nmid h$.

\begin{proof}[Proof of Lemma~\ref{lem:divisibility}]
We closely follow the proof of~\cite[Proposition~5.1]{Bloom2025density}. The only new ingredient is Lemma~\ref{lem:shiftedsieve}, which uses that $A$ consists of shifted primes: it shows that the sets $A^{(1)}$ and $A^{(2)}$ defined below have small reciprocal sums, whereas Bloom in~\cite{Bloom2025density} instead passes to one of these sets. Throughout the proof, we assume that $N$ is sufficiently large in terms of $h$.

Fix $t\in\mathfrak{m}_{2}$. Since $t$ is fixed, we write $I=I_{t}$ and $A_{I}=A_{I_{t}}$, so that $A_{I}$ is the set of those $n\in A$ that divide some element of $I$. As $t\notin\mathfrak{m}_{1}$, we have $|A\setminus A_{I}|<M/\log N$. Moreover, $A_{q}\setminus A_{I}=\{n\in A_{q}:\ \|t/n\|\geq K/(2n)\}$ for every $q\in\mathcal{Q}_{A}$, so by the definition of $\mathcal{D}_{t}$,
\[
|A_{q}\setminus A_{I}|<L\qquad\text{for all } q\in\mathcal{D}_{t}.
\]

Let $\mathcal{E}_{I}$ be the set of those $q\in\mathcal{Q}_{A}$ such that $R(A_{I};q)>(\log N)^{-1/100}/2$, where \[R(A_{I};q)=\sum_{n\in A_{q}\cap A_{I}}\frac{q}{n}.\]
For every $q\in\mathcal{D}_{t}$, every element of $A_{q}\setminus A_{I}$ is greater than $M$, so by condition (iii) of Proposition~\ref{prop:main},
\[
R(A_{I};q)=R(A;q)-\sum_{n\in A_{q}\setminus A_{I}}\frac{q}{n}\geq \eta-\frac{qL}{M}.
\]
Since $q\leq N^{0.99}$, $\tau(1-\tau)\geq\frac{1}{2}(\log N)^{-2}$, $K=MN^{-2/\log_{2}N}$ and $M\geq N^{0.997}$, we have
\[
\frac{qL}{M}=\frac{25\,qN^{2}\log N\log_{2}N}{\tau(1-\tau)K^{2}M}\leq\frac{50\,N^{0.99}N^{2+4/\log_{2}N}(\log N)^{3}\log_{2}N}{M^{3}}\leq N^{-0.001+o(1)}\leq\frac{\eta}{2}.
\]
Hence $R(A_{I};q)\geq \eta/2$, and so in particular $\mathcal{D}_{t}\subseteq\mathcal{E}_{I}$.

For $q\in\mathcal{E}_{I}$ we have $R(A_{I};q)\geq(\log N)^{-1/2}$, so we may apply Lemma~\ref{lem-usingq} with $A$ replaced by $A_{I}$ and $k=\log_{2}N/\log(2\log_{2}N)$. Since $(\log N)^{1-1/k}={\log N}/{(2\log_{2}N)}$, this produces some $d_{q}$ such that
\[
qd_{q}>M\exp\Big(-\frac{\log N}{2\log_{2}N}\Big)=MN^{-1/(2\log_{2}N)}>K
\]
and $\omega(d_{q})<\frac{1}{500}\log_{2}N$, and
\[
\sum_{\substack{n\in A_{q}\cap A_{I},\ qd_{q}\mid n\\ (qd_{q},n/qd_{q})=1}}\frac{qd_{q}}{n}\gg\frac{R(A_{I};q)}{(2\log_{2}N)^{2}}\gg\frac{1}{(\log N)^{1/100}(\log_{2}N)^{2}}.
\]
By the definition of $A_{I}$, every $n\in A_{I}$ with $qd_{q}\mid n$ must divide some element of $I$. In fact, they all divide the same element of $I$, which we call $x_{q}$: each such $n$ is divisible by $qd_{q}>K$, and an integer greater than $K$ divides at most one element of $I$. In particular, $q\mid x_{q}$.

Let
\[
A_{I}^{(q)}=\Big\{\frac{n}{qd_{q}}:\ n\in A_{q}\cap A_{I}\text{ with } qd_{q}\mid n\text{ and }\Big(qd_{q},\frac{n}{qd_{q}}\Big)=1\Big\},
\]
so that, assuming $N$ is sufficiently large, $R(A_{I}^{(q)})\geq(\log N)^{-1/99}$, say. Each $m=n/(qd_{q})\in A_{I}^{(q)}$ satisfies $\omega(m)\geq\omega(n)-1-\omega(d_{q})\geq\big(\frac{99}{100}-\frac{1}{500}\big)\log_{2}N-1\geq\frac{97}{99}\log_{2}N$ and $\omega(m)\leq 2\log_{2}N$. So we may apply Lemma~\ref{lem-rtop} with $\varepsilon=\frac{2}{99}$, which gives
\[
\sum_{r\in\mathcal{Q}_{A_{I}^{(q)}}}\frac{1}{r}\geq\frac{95}{99}e^{-1}\log_{2}N.
\]
Trivially, $\mathcal{Q}_{A_{I}^{(q)}}\subseteq\mathcal{Q}_{A}$. Moreover, by the choice of $x_{q}$, every $r\in\mathcal{Q}_{A_{I}^{(q)}}$ divides $x_{q}$, and hence
\begin{equation}\label{eq:xq}
\sum_{\substack{r\mid x_{q}\\ r\in\mathcal{Q}_{A}}}\frac{1}{r}\geq\frac{95}{99}e^{-1}\log_{2}N\geq 0.35\log_{2}N.
\end{equation}
For any two $n_{1}\neq n_{2}\in I$, we have $0<|n_{1}-n_{2}|<K\leq N$, and hence
\begin{equation}\label{eq:overlap}
\sum_{q\mid(n_{1},n_{2})}\frac{1}{q}\leq C\log_{3}N\leq 0.01\log_{2}N
\end{equation}
for $N$ sufficiently large, by Lemma~\ref{lem-basic}; here $C$ is the implied constant in that lemma.

Since every $q\in\mathcal{Q}_{A}$ is a prime power at most $N$, Mertens' estimate \eqref{mertens1} gives
\[
\sum_{q\in\mathcal{Q}_{A}}\frac{1}{q}\leq(1+o(1))\log_{2}N.
\]
It follows that there are at most two distinct values of $x_{q}\in I$ as $q$ ranges over $\mathcal{E}_{I}$. Indeed, suppose that there were three distinct such values $w_{1}$, $w_{2}$ and $w_{3}$. Then, by \eqref{eq:xq} and \eqref{eq:overlap},
\begin{align*}
(1+o(1))\log_{2}N\geq\sum_{q\in\mathcal{Q}_{A}}\frac{1}{q}&\geq\sum_{i=1}^{3}\ \sum_{\substack{q\mid w_{i}\\ q\in\mathcal{Q}_{A}}}\frac{1}{q}-\sum_{1\leq i<j\leq 3}\ \sum_{q\mid(w_{i},w_{j})}\frac{1}{q}\\
&\geq 3(0.35\log_{2}N)-3(0.01\log_{2}N)=1.02\log_{2}N,
\end{align*}
which is a contradiction.

If all the $x_{q}$ with $q\in\mathcal{E}_{I}$ are equal, we take $x(t)$ to be their common value; then every $q\in\mathcal{D}_{t}\subseteq\mathcal{E}_{I}$ divides $x(t)$. (If $\mathcal{E}_{I}$ is empty, any $x(t)\in I$ works.) It remains to rule out the case of exactly two values, say $w_{1}=x_{q_{1}}$ and $w_{2}=x_{q_{2}}$ with $q_{1},q_{2}\in\mathcal{E}_{I}$ and $w_{1}\neq w_{2}$. We assume that this case occurs and derive a contradiction. Let
\[
A^{(i)}=\{n\in A:\ n\mid w_{i}\}\quad(i=1,2)\qquad\text{and}\qquad A^{(0)}=A\setminus(A^{(1)}\cup A^{(2)}).
\]

We first show that $A^{(1)}$ and $A^{(2)}$ have small reciprocal sums; this is where we use that $A$ consists of shifted primes. Let $\mathcal{P}^{(1)}$ be the set of primes $p$ with $3\leq p\leq M^{1/10}$, $p\nmid h$ and $p\nmid w_{2}$ such that $p\mid r$ for some $r\in\mathcal{Q}_{A}$ with $r\mid w_{1}$. To bound $\sum_{p\in\mathcal{P}^{(1)}}1/p$ from below, we start from \eqref{eq:xq} with $q=q_{1}$. Prime powers with exponent at least $2$ contribute $O(1)$, and the primes $p=2$, $p\mid h$ and $M^{1/10}<p\leq N$ contribute $O_{h}(1)$, the last range by Mertens' estimate since $\log M\geq 0.997\log N$. The primes dividing both $w_{1}$ and $w_{2}$ contribute at most $\sum_{q\mid(w_{1},w_{2})}1/q$, which is $O(\log_{3}N)$ by \eqref{eq:overlap}. Hence
\[
\sum_{p\in\mathcal{P}^{(1)}}\frac{1}{p}\geq 0.35\log_{2}N-O_{h}(1)-\sum_{q\mid(w_{1},w_{2})}\frac{1}{q}\geq 0.35\log_{2}N-O_{h}(1)-O(\log_{3}N)\geq 0.34\log_{2}N.
\]
Every $n\in A^{(2)}$ divides $w_{2}$, and hence $p\nmid n$ for all $p\in\mathcal{P}^{(1)}$. Moreover, $n+h$ is prime for every $n\in A$.

Let $j_{0}=\lfloor\log(N/M)/\log 2\rfloor$ and $X_{j}=2^{j}M$ for $0\leq j\leq j_{0}$, so that $(M,N]\subseteq\bigcup_{j=0}^{j_{0}}[X_{j},2X_{j})$. For each $j$, $X_{j}\geq M$, and every $p\in\mathcal{P}^{(1)}$ satisfies $3\leq p\leq M^{1/10}\leq X_{j}^{1/10}$ and $p\nmid h$. Hence, by Lemma~\ref{lem:shiftedsieve},
\begin{align*}
|A^{(2)}\cap[X_{j},2X_{j})|&\leq C_{2}\frac{X_{j}}{\log X_{j}}\prod_{p\in\mathcal{P}^{(1)}}\Big(1-\frac{1}{p-1}\Big)\\
&\leq C_{2}\frac{X_{j}}{\log M}\exp\Big(-\sum_{p\in\mathcal{P}^{(1)}}\frac{1}{p}\Big)\leq C_{2}\frac{X_{j}}{\log M}(\log N)^{-0.34}.
\end{align*}
Since $j_{0}\ll\log N$ and $\log M\geq 0.997\log N$, it follows that
\[
R(A^{(2)})\leq\sum_{j=0}^{j_{0}}\frac{|A^{(2)}\cap[X_{j},2X_{j})|}{X_{j}}\leq\frac{(j_{0}+1)C_{2}}{\log M}(\log N)^{-0.34}\ll(\log N)^{-0.34}.
\]
Similarly, $R(A^{(1)})\ll(\log N)^{-0.34}$. Since $R(A)\geq(\log N)^{-1/101}$, we conclude that, for $N$ sufficiently large,
\[
R(A^{(0)})\geq\frac{R(A)}{2}.
\]

Let $A'$ be the set of $n\in A_{I}\cap A^{(0)}$ such that every $q\in\mathcal{Q}_{A}$ with $q\pdiv n$ lies in $\mathcal{E}_{I}$. Every $n\in A^{(0)}\setminus A'$ either lies in $A\setminus A_{I}$ or lies in $A_{q}\cap A_{I}$ for some $q\in\mathcal{Q}_{A}\setminus\mathcal{E}_{I}$. Hence, by the definition of $\mathcal{E}_{I}$, the bound $|A\setminus A_{I}|<M/\log N$ and Mertens' estimate \eqref{mertens1},
\[
R(A^{(0)}\setminus A')\leq\frac{|A\setminus A_{I}|}{M}+\sum_{q\in\mathcal{Q}_{A}\setminus\mathcal{E}_{I}}\frac{1}{q}R(A_{I};q)\ll\frac{\log_{2}N}{(\log N)^{1/100}}.
\]
Since $R(A^{(0)})\geq R(A)/2\geq(\log N)^{-1/101}/2$, it follows that $R(A')\gg(\log N)^{-1/101}$.

In particular, since every element of $A'$ is greater than $M$,
\[
|A'|\geq MR(A')\gg\frac{M}{(\log N)^{1/101}}.
\]
Every $n\in A'\subseteq A_{I}$ divides some element of $I$, and $|I|\leq K+1\leq 2K$. Therefore, there must exist some $v\in I$ (necessarily $v\neq w_{1}$ and $v\neq w_{2}$ since $A'\subseteq A^{(0)}$) such that, if $A''=\{n\in A':\ n\mid v\}$, then
\[
|A''|\geq\frac{|A'|}{2K}\gg N^{2/\log_{2}N}(\log N)^{-1/101},
\]
and hence $|A''|\geq N^{3/(2\log_{2}N)}$, say.

However, if $n\in A''$, then both $n\mid v$ and $n\mid w_{1}w_{2}$ (since every $q$ with $n\in A_{q}$ is in $\mathcal{E}_{I}$ and so divides either $w_{1}$ or $w_{2}$), and hence $n$ divides
\[
(v,w_{1}w_{2})\leq(v,w_{1})(v,w_{2})\leq|v-w_{1}||v-w_{2}|\leq N^{2}.
\]
Therefore, the size of $A''$ is at most the number of divisors of the fixed integer $(v,w_{1}w_{2})\leq N^{2}$, which is at most $N^{(1+o(1))2\log 2/\log_{2}N}$ (see~\cite[Theorem~2.11]{MV}), and hence we have a contradiction for large enough $N$, since $2\log 2<3/2$.   Thus the $x_{q}$ with $q\in\mathcal{E}_{I}$ cannot take two distinct values. So they are all equal, which completes the proof of Lemma~\ref{lem:divisibility}.
\end{proof}

\bibliographystyle{plain}
\bibliography{Golombprimes}
\end{document}